    \documentstyle [righttag] {amsart}

    \input amssym.def
    \input amssym

    \newtheorem{propo}{Proposition}
    \newtheorem{theo}[propo]{Theorem}
    \newtheorem{lemma}[propo]{Lemma}
    \newtheorem{corol}[propo]{Corollary}

    \theoremstyle{defi}
    \newtheorem{defi}[propo]{Definition}

    \theoremstyle{rema}
    \newtheorem{rema}[propo]{Remark}

\numberwithin{equation}{section}
\numberwithin{propo}{section}

\begin{document}
    \title[Local systems of twisted vertex operators]
{Local systems of twisted vertex operators,
vertex operator superalgebras and twisted modules}
    \author{Hai-sheng Li}
    \address{Department of Mathematics\\ Rutgers University
\\ New Brunswick, NJ 08903}
    \curraddr{Department of Mathematics\\University of California\\
Santa Cruz, CA 95064}
    \email{hli@@math.ucsc.edu}
   \thanks{This paper is in final form and no version of it will
be submitted for publication elsewhere.}
	\bibliographystyle{alpha}
	\maketitle

	\newcommand{\nno}{\nonumber}
	\newcommand{\lbar}{\bigg\vert}
	\newcommand{\p}{\partial}
	\newcommand{\dps}{\displaystyle}
	\newcommand{\bra}{\langle}
	\newcommand{\ket}{\rangle}
 \newcommand{\res}{\mbox{\normalshape Res}}
 \newcommand{\epf}{\hspace{2em}$\Box$}
 \newcommand{\epfv}{\hspace{1em}$\Box$\vspace{1em}}
\newcommand{\nord}{\mbox{\scriptsize ${\circ\atop\circ}$}}
\newcommand{\wt}{\mbox{\normalshape wt}\ }

\begin{abstract}
We introduce the
notion of ``local system of $\Bbb{Z}_{T}$-twisted vertex operators'' on
a $\Bbb{Z}_{2}$-graded vector space $M$, generalizing the
notion of local system of vertex operators [Li]. First, we prove
that any local system of $\Bbb{Z}_{T}$-twisted vertex operators
on $M$ has a vertex superalgebra structure with an
automorphism $\sigma$ of order $T$ with $M$ as a $\sigma$-twisted
module. Then we prove that
for a vertex (operator) superalgebra $V$ with an automorphism $\sigma$ of order
$T$, giving a $\sigma$-twisted $V$-module $M$ is equivalent to giving
a vertex (operator) superalgebra homomorphism from $V$ to some local
system of $\Bbb{Z}_{T}$-twisted vertex operators on $M$. As
applications, we study the twisted modules for vertex operator
(super)algebras associated to some well-known infinite-dimensional Lie
(super)algebras and we prove the complete reducibility of
$\Bbb{Z}_{T}$-twisted modules
for vertex operator algebras associated to standard modules for
an affine Lie algebra.
\end{abstract}

\section{Introduction}
In [Li], motivated by [G] we introduced the notion of local system of vertex
operators
and proved that any local system of vertex operators on a
super vector space $M$ has a vertex
superalgebra structure with $M$ as a module, so that if $V$ is
a vertex (operator) superalgebra, then giving a $V$-module $M$ is
equivalent to giving a vertex superalgebra homomorphism from $V$ to
some local system of vertex operators on $M$. Similar to the notion of local
system of
vertex operators, the notion of commutative quantum operator algebra
was introduced in [LZ1-2]. For a  vertex (operator)
superalgebra $V$, in addition to having module theory, more generally
we have twisted module theory ([D], [FFR], [FLM]). In this paper,
generalizing the notion of local system of vertex operators ([G], [Li]), we
introduce the
notion of what we call ``local system of $\Bbb{Z}_{T}$-twisted vertex
operators'' on a $\Bbb{Z}_{2}$-graded vector space $M$ and we prove similar
results.

Let  $M$ be a $\Bbb{Z}_{2}$-graded vector space and
let $T$ be a fixed positive integer.
A {\it weak $\Bbb{Z}_{T}$-twisted vertex operator} on $M$ is a formal
series $a(z)\in ({\rm End}M)[[z^{{1\over T}},z^{-{1\over T}}]]$
satisfying $a(z)u\in M((z^{{1\over T}}))$ for any $u\in M$. Let
$F(M,T)$ be the vector space of all weak $\Bbb{Z}_{T}$-twisted
vertex operators on $M$.
Set $\varepsilon=\exp
\left(\frac{2\pi \sqrt{-1}}{T}\right)$. Let $\sigma$ be the endomorphism of
$({\rm End}M)[[z^{{1\over T}},z^{-{1\over T}}]]$ defined by: $\sigma
f(z^{1\over T})=f(\varepsilon^{-1} z^{1\over T})$. Then $\sigma$ is an
automorphism of $F(M,T)$ of order $T$. For any integer $n$, set
\begin{eqnarray}
F(M,T)^{n}=\{f(z^{{1\over T}})\in F(M,T)|\sigma f(z^{{1\over
T}})=\varepsilon^{n}f(z^{{1\over T}})\}.
\end{eqnarray}
Then $z^{{n\over T}}F(M,T)^{n}\subseteq  ({\rm End}M)[[z,z^{-1}]]$
 and
$F(M,T)^{n}=F(M,T)^{T+n}$ for any $n\in \Bbb{Z}$. Thus
\begin{eqnarray}
F(M,T)=F(M,T)^{0}\oplus F(M,T)^{1}\oplus \cdots \oplus F(M,T)^{T-1},
\end{eqnarray}
so that $F(M,T)$ is a $\left(\Bbb{Z}_{2} \times \Bbb{Z}_{T}\right)$-graded
vector space.

We define {\it a local subspace} of $F(M,T)$ to be a
$\left(\Bbb{Z}_{2}\times \Bbb{Z}_{T}\right)$-graded subspace $A$
such that for any two $\Bbb{Z}_{2}$-homogeneous elements $a(z)$ and $b(z)$ of
$A$, there is a positive integer $m$ such that
\begin{eqnarray}
(z_{1}-z_{2})^{m}a(z_{1})b(z_{2})=(-1)^{|a(z)||b(z)|}(z_{1}-z_{2})^{m}
b(z_{2})a(z_{1}),
\end{eqnarray}
where $|a(z)|=0$ if $a(z)$ is even, $|a(z)|=1$ if $a(z)$ is odd.
And we define {\it a local system of $\Bbb{Z}_{T}$-twisted
vertex operators on $M$} to be a maximal local (graded)
subspace of $F(M,T)$.

Let $M$ be a $\Bbb{Z}_{2}$-graded vector space and let $A$ be
a local system of $\Bbb{Z}_{T}$-twisted vertex operators on
$M$. Then $\sigma$ is an automorphism of $A$ such that
$\sigma^{T}=Id_{A}$. Set $A^{k}=A\cap F(M,T)^{k}$ for any
integer $k$. Similar to the untwisted case, the ``multiplication'' for
twisted vertex operators comes from the
twisted iterate formula. For any $a(z)\in
A^{k}$, $b(z)\in A$ and any integer $n\in \Bbb{Z}$,
$a(z)_{n}b(z)$ as an element of $F(M,T)$ is defined by:
\begin{eqnarray}
& &\;\;\;Y(a(z),z_{0})b(z)\nonumber\\
& &=:\sum_{n\in \Bbb{Z}}(a(z)_{n}b(z))z_{0}^{-n-1}\nonumber\\
& &={\rm Res}_{z_{1}}\!\left(\!\frac{z_{1}-z_{0}}{z}\!\right)^{{k\over T}}
\!\left(\!z_{0}^{-1}\delta\!\left(\!\frac{z_{1}-z}{z_{0}}\!\right)\!a(z_{1})b(z)
-z_{0}^{-1}\delta\!\left(\!\frac{z-z_{1}}{-z_{0}}\!\right)\!b(z)a(z_{1})\right).
\nonumber\\
& &\mbox{}
\end{eqnarray}
As in the untwisted case [Li], we first prove that $A$ is closed
under the  ``multiplication'' (1.4) (Proposition 3.9). Then we prove
that $A$ under the ``adjoint action'' (1.4) is a local system of
(untwisted) vertex operators on $A$ (Proposition 3.13). By Proposition 2.5
(recalled from [Li]), $A$ is a vertex superalgebra with a natural automorphism
$\sigma$ such that $\sigma^{T}=Id_{A}$. Furthermore, since the
supercommutativity and the twisted iterate formula imply the twisted Jacobi
identity (for a twisted module), $M$
is a $\sigma$-twisted module for the vertex superalgebra (Theorem 3.14).
Let $V$ be a vertex superalgebra with an automorphism $\sigma$ of order $T$.
Then we prove that giving a $\sigma$-twisted $V$-module $M$ is equivalent
to giving a vertex superalgebra homomorphism from $V$ to some local
system of $\Bbb{Z}_{T}$-twisted vertex operators on $M$
(Proposition 3.17).

Let $M$ be any $\Bbb{Z}_{2}$-graded vector space and let $S$ be
any set of mutually local homogeneous $\Bbb{Z}_{T}$-twisted vertex
operators on $M$. It follows from Zorn's lemma that there exists a
local system $A$ containing $S$. Therefore,
$S$ generates (inside $A$ under the twisted ``multiplication'' (1.4)) a vertex
superalgebra with $\sigma$ as an automorphism
such that $M$ is a
$\sigma$-twisted module (Corollary 3.15).  Since the
``multiplication'' (1.4) does not depend on the choice of the local
system $A$, the vertex superalgebra $\langle S\rangle$ is canonical,
so that we can speak about the vertex superalgebra generated by $S$.

To describe our results, for simplicity, we consider a special case.
Let ${\bf g}$ be a finite-dimensional simple Lie algebra with a fixed
Cartan subalgebra ${\bf h}$. Let $\langle\cdot,\cdot\rangle$ be the normalized
Killing form on ${\bf g}$ such that the square of the length of the
longest root is $2$. It has been proved ([DL], [FZ], [Li], [Lia]) that the
generalized Verma $\tilde{{\bf g}}$-module $M_{{\bf g}}(\ell,{\bf C})$
has a natural vertex operator algebra structure for any $\ell\ne -\Omega$.
Let $\sigma$ be an automorphism of ${\bf g}$ of
order $T$. Let $M$ be any restricted module of level $\ell$
for the twisted affine Lie algebra $\tilde{{\bf g}}[\sigma]$.
By considering the generating function $a(z)=\sum_{n\in \Bbb{Z}}(t^{n+{j\over
T}})z^{-n-1-{j\over T}}$ for $a\in {\bf
g}^{j}$ as an element of $F(M,T)$, we obtain a local
subspace  $\{ a(z)|a\in {\bf g}\}$, so that $\{ a(z)|a\in {\bf g}\}$
generates a vertex
algebra $V$ with $M$ as a $\sigma$-twisted module. Since
$V$ is a quotient vertex algebra of $M_{{\bf g}}(\ell,{\bf
C})$, $M$ is a weak $\sigma$-twisted  $M_{{\bf g}}(\ell,{\bf
C})$-module.

Let ${\bf g}$ be a finite-dimensional simple Lie algebra and let
$\ell$ be a positive integer. It has been proved ([DL], [FZ],
[Li]) that any lower truncated $\Bbb{Z}$-graded weak $L_{{\bf
g}}(\ell,0)$-module is a direct sum of standard $\tilde{{\bf
g}}$-modules of level $\ell$ and that the set of equivalence classes of
irreducible $L_{{\bf g}}(\ell,0)$-modules is exactly the set of
equivalence classes of standard $\tilde{{\bf g}}$-modules of level $\ell$.
Let $\mu$ be a Dynkin diagram automorphism of ${\bf g}$. Then
using a similar method used for ordinary $V$-modules ([DL], [Li]) we prove that
any lower truncated ${1\over T}\Bbb{Z}$-graded
weak $\mu$-twisted $L_{{\bf g}}(\ell,0)$-module is
completely reducible and that the set of irreducible $\mu$-twisted $L_{{\bf
g}}(\ell,0)$-modules
is exactly the set of the standard $\tilde{{\bf g}}[\mu]$-modules  of
level $\ell$ (Theorem 5.8).

Let $\sigma$ be any automorphism of order $T$ of ${\bf g}$. By [K],
$\sigma$ is conjugate to an automorphism $\mu \sigma_{h}$, where $\mu$
is a Dynkin diagram automorphism of ${\bf g}$ and $\sigma_{h}$ is an
inner automorphism of ${\bf g}$ associated with an element $h$ of the
Cartan subalgebra ${\bf h}$ such that $\sigma (h)=h$. If two
automorphisms $\sigma_{1}$ and $\sigma_{2}$ of a vertex operator algebra $V$
are conjugate each other, it is clear that the equivalence classes of
$\sigma_{1}$-twisted $V$-modules one-to-one correspond to the
equivalence classes of $\sigma_{2}$-twisted $V$-modules.
In the Lie algebra level, it is well-known ([K], [H]) that the
twisted affine Lie algebras $\tilde{{\bf g}}[\mu]$ and $\tilde{{\bf
g}}[\mu\sigma_{h}]$ are isomorphic each other, so that the equivalence
classes of
$\tilde{{\bf g}}[\mu]$-modules one-to-one correspond to the
equivalence classes of
$\tilde{{\bf g}}[\mu\sigma_{h}]$-modules. But it is not obvious
that we have a one-to-one correspondence between the equivalence
classes of $\mu$-twisted $L_{{\bf g}}(\ell,0)$-modules and the
equivalence classes of $(\mu\sigma_{h})$-twisted
 $L_{{\bf g}}(\ell,0)$-modules.
In Section 5, such a one-to-one correspondence is established as a
corollary of Proposition 5.4. Moreover, Proposition 5.4 has also its
own interest. We will study certain applications of Proposition 5.4 in
a coming paper.
In this way, we show that for any automorphism $\sigma$ of ${\bf g}$,
any lower truncated $({1\over T}\Bbb{Z})$-graded weak $\sigma$-twisted
$L_{{\bf g}}(\ell,0)$-module is completely reducible and that there are
only finitely many irreducible $\sigma$-twisted $L_{{\bf
g}}(\ell,0)$-modules up to equivalence.

This paper is organized as follows: Section 2 is a preliminary section.
In Section 3 we introduce local systems of
$\Bbb{Z}_{T}$-twisted vertex operators.
In Section 4 we study $\Bbb{Z}_{T}$-twisted modules for
vertex operator superalgebras  associated to the Neveu-Schwarz algebra
and an affine Lie superalgebra. In Section 5
we give the semisimple twisted module theory
for vertex operator algebras associated to standard
modules for an affine Lie algebra.

{\bf Acknowledgment}
We would like to thank Professors Chongying Dong, James Lepowsky and
Robert Wilson for many useful discussions.

\section{Vertex superalgebras and twisted modules}
In this section, we review the definitions of vertex (operator)
superalgebra and twisted module and we also present some elementary results.

Let $z$, $z_{0}$, $z_{1},\cdots$ be
commuting formal variables. For a vector space $V$, we recall the
following notations and formal variable calculus from [FLM]:
\begin{eqnarray}
& &V\{z\}=\{\sum_{n\in {\bf C}}v_{n}z^{n}|v_{n}\in V\},\;
V[[z,z^{-1}]]=\{\sum_{n\in \Bbb{Z}}v_{n}z^{n}|v_{n}\in V\},\\
& &V[z]=\{\sum_{n\in {\bf N}}v_{n}z^{n}|v_{n}\in V,v_{n}=0\;\mbox{ for }
n\mbox{ sufficiently large}\},\\
& &V[z,z^{-1}]=\{\sum_{n\in \Bbb{Z}}v_{n}z^{n}|v_{n}\in V,v_{n}=0\;\mbox{ for
all but finitely many }n\},\\
& &V[[z]]=\{\sum_{n\in{\bf N}}v_{n}z^{n}|v_{n}\in V\},\\
& &V((z))=\{\sum_{n\in \Bbb{Z}}v_{n}z^{n}|v_{n}\in V,v_{n}=0\;\mbox{
for }n\;\mbox{sufficiently small}\}.
\end{eqnarray}
For $f(z)=\sum_{n\in\Bbb{Z}}v_{n}z^{n}\in V[[z,z^{-1}]]$, the formal
derivative is defined to be
\begin{eqnarray}
{d\over dz}f(z)=f'(z)=\sum_{n\in \Bbb{Z}}nv_{n}z^{n-1},
\end{eqnarray}
and the formal residue is defined as follows:
\begin{eqnarray}
{\rm Res}_{z}f(z)=v_{-1}\;\;\;\mbox{(the coefficient of $z^{-1}$ in $f(z)$)}.
\end{eqnarray}
If $f(z)\in {\bf C}((z)), g(z)\in V((z))$, we have:
\begin{eqnarray}
{\rm Res}_{z}(f'(z)g(z))=-{\rm Res}_{z}(f(z)g'(z)).
\end{eqnarray}
For $\alpha\in \Bbb{C}$, as a formal power series,
$(z_{1}+z_{2})^{\alpha}$ is defined to be
\begin{eqnarray}
(z_{1}+z_{2})^{\alpha}=\sum_{k=0}^{\infty}\left(\begin{array}{c}\alpha\\k
\end{array}\right)z_{1}^{\alpha -k}z_{2}^{k},
\end{eqnarray}
where
$\displaystyle{\left(\begin{array}{c}\alpha\\k\end{array}\right)=\frac{\alpha
(\alpha-1)\cdots (\alpha-k+1)}{k!}}$.
If $f(z)\in V\{z\}$, then we have the following Taylor formula:
\begin{eqnarray}
e^{z_{0}{\partial\over\partial z}}f(z)=f(z+z_{0}).
\end{eqnarray}
The formal $\delta$-function is defined to be
\begin{eqnarray}
\delta(z)=\sum_{n\in \Bbb{Z}}z^{n}\in \Bbb{C}[[z,z^{-1}]].
\end{eqnarray}
Thus
\begin{eqnarray}
\delta\left(\frac{z_{1}-z_{2}}{z_{0}}\right)=\sum_{n\in \Bbb{Z}}z_{0}^{-n}
(z_{1}-z_{2})^{n}=\sum_{n\in \Bbb{Z}}\sum_{k\in \Bbb{N}}(-1)^{k}\left(
\begin{array}{c}n\\k\end{array}\right)z_{0}^{-n}z_{1}^{n-k}z_{2}^{k}.
\end{eqnarray}
Furthermore, for $\alpha\in \Bbb{C}$, we have:
\begin{eqnarray}
z_{0}^{-1}\delta\left(\frac{z_{1}-z_{2}}{z_{0}}\right)\left(\frac{z_{1}
-z_{2}}{z_{0}}\right)^{\alpha}=z_{1}^{-1}\delta\left(\frac{z_{0}+z_{2}}{z_{1}}
\right)\left(\frac{z_{0}+z_{2}}{z_{1}}\right)^{-\alpha}.
\end{eqnarray}

\begin{lemma}\label{fhl}
 [FLM] If $f(z_{1},z_{2})\in
V[[z_{1},z_{1}^{-1}]]((z_{2}))$, then
\begin{eqnarray}
\delta\left(\frac{z_{0}+z_{2}}{z_{1}}\right)f(z_{1},z_{2})=
\delta\left(\frac{z_{0}+z_{2}}{z_{1}}\right)f(z_{0}+z_{2},z_{2}).
\end{eqnarray}
\end{lemma}

\begin{lemma}
Let $m,n\in \Bbb{Z}_{+}$ such that $m>n$. Then
\begin{eqnarray}
(z_{1}-z_{2})^{m}\delta ^{(n)}\left({z_{1}\over z_{2}}\right)=0.
\end{eqnarray}
\end{lemma}

{\bf Proof.} It easily follows from Lemma 2.1 and an induction
on $n$.$\;\;\;\;\Box$

\begin{lemma}
 Let $V$ be any vector space, let $\alpha$ and $\beta$ be
rational numbers and let $f_{j}(z_{2})$
($j=0,\cdots,n$) be elements of $V((z_{2}^{\beta}))$. Let
\begin{eqnarray}
g(z_{1},z_{2})=z_{2}^{-1}\delta\left({z_{1}\over
z_{2}}\right)\left({z_{1}\over z_{2}}\right)^{\alpha}.
\end{eqnarray}
Then
\begin{eqnarray}
g(z_{1},z_{2})f_{0}(z_{2})\!+\!\left({\partial\over\partial
z_{2}}g(z_{1},z_{2})
\right)\!f_{1}(z_{2})\!+\!\cdots\!+\!\left({\partial^{n}\over\partial
z_{2}^{n}}
g(z_{1},z_{2})\right)\!f_{n}(z_{2})=0
\end{eqnarray}
if and only if $f_{j}(z_{2})=0$ for all $j$.
\end{lemma}

{\bf Proof.} First, we prove that for any nonnegative integer $k$ we have:
\begin{eqnarray}
\left({\partial \over\partial z_{2}}\right)^{k}g(z_{1},z_{2})
=(-1)^{k}\left({\partial\over\partial z_{1}}\right)^{k}g(z_{1},z_{2}).
\end{eqnarray}
Using definition, we have:
\begin{eqnarray}
{\partial\over\partial z_{2}}g(z_{1},z_{2})
=-{\partial\over\partial z_{1}}g(z_{1},z_{2}).
\end{eqnarray}
Then (2.18) easily follows from an induction on $k$.
Suppose $f_{j}(z_{2})\ne 0$ for some $j$. Without losing generality,
we may assume $f_{n}(z_{2})\ne 0$.
Applying ${\rm Res}_{z_{1}}z_{1}^{-\alpha}(z_{1}-z_{2})^{n}$ to (2.17),
then using Leibniz rule and Lemma 2.2 we obtain:
\begin{eqnarray}
0&=&{\rm Res}_{z_{1}}z_{1}^{-\alpha}(z_{1}-z_{2})^{n}
\left({\partial^{n}\over\partial z_{2}^{n}}
g(z_{1},z_{2})\right)f_{n}(z_{2})\nonumber\\
&=&(-1)^{n}{\rm Res}_{z_{1}}z_{1}^{-\alpha}(z_{1}-z_{2})^{n}
\left({\partial^{n}\over\partial z_{1}^{n}}
g(z_{1},z_{2})\right)f_{n}(z_{2})\nonumber\\
&=&{\rm Res}_{z_{1}}\left({\partial^{n}\over\partial z_{1}^{n}}
z_{1}^{-\alpha}(z_{1}-z_{2})^{n}\right)g(z_{1},z_{2})f_{n}(z_{2})\nonumber\\
&=&n!z_{2}^{-\alpha}f_{n}(z_{2}).
\end{eqnarray}
It is a contradiction. $\;\;\;\;\Box$

Let $M=M^{(0)}\oplus M^{(1)}$ be any $\Bbb{Z}_{2}$-graded vector
space. Then any element $u$ in $M^{(0)}$ (resp. $M^{(1)}$) is said to be
{\it even} (resp. {\it odd}).
For any homogeneous element $u$, we define $|u|=0$ if $u$
is even, $|u|=1$ if $u$ is odd. If $M$ and $W$ are any two
$\Bbb{Z}_{2}$-graded vector spaces, we define
$\varepsilon_{u,v}=(-1)^{|u||v|}$ for any homogeneous elements $u\in
M, v\in W$.

\begin{defi}
 A {\it vertex superalgebra} is a quadruple
$(V,{\bf 1},D,Y)$, where $V=V^{(0)}\oplus V^{(1)}$ is a $\Bbb{Z}_{2}$-graded
vector space, $D$ is a $\Bbb{Z}_{2}$-endomorphism
of $V$, ${\bf 1}$ is a specified vector called the {\it vacuum} of $V$,
and $Y$ is a linear map
\begin{eqnarray}
Y(\cdot,z):& &V\rightarrow ({\rm End}V)[[z,z^{-1}]];\nonumber\\
& &a\mapsto Y(a,z)=\sum_{n\in \Bbb{Z}}a_{n}z^{-n-1}\;\;(\mbox{where }
a_{n}\in {\rm End}V)
\end{eqnarray}
such that
\begin{eqnarray*}
(V1)& &\mbox{For any }a,b\in V, a_{n}b=0\;\;\;\mbox{ for }n
\mbox{ sufficiently large;}\\
(V2)& &[D,Y(a,z)]=Y(D(a),z)={d\over dz}Y(a,z)\;\;\mbox{  for any }a\in V;\\
(V3)& &Y({\bf 1},z)=Id_{V}\;\;\;\mbox{(the identity operator of $V$)};\\
(V4)& &Y(a,z){\bf 1}\in ({\rm End}V)[[z]] \mbox{ and }\lim_{z \rightarrow
0}Y(a,z){\bf 1}=a\;\;\mbox{  for any }a\in V;\\
(V5)& &\mbox{For }\Bbb{Z}_{2}\mbox{ -homogeneous }a,b\in V,
\mbox{ the following {\it Jacobi identity} holds:}
\end{eqnarray*}
\begin{eqnarray}
&
&\;\;\;z_{0}^{-1}\delta\left(\frac{z_{1}-z_{2}}{z_{0}}\right)Y(a,z_{1})Y(b,z_{2})
-\varepsilon_{a,b}z_{0}^{-1}\delta\left(\frac{z_{2}-z_{1}}{-z_{0}}\right)
Y(b,z_{2})Y(a,z_{1})\nonumber \\
& &=z_{2}^{-1}\delta\left(\frac{z_{1}-z_{0}}{z_{2}}\right)Y(Y(a,z_{0})b,z_{2}).
\end{eqnarray}
\end{defi}

This completes the definition of vertex superalgebra.
 For any $\Bbb{Z}_{2}$-homogeneous elements $a, b\in V$, there is a
positive integer $m$ such that (cf. [DL], [Li])
\begin{eqnarray}
(z_{1}-z_{2})^{m}Y(a,z_{1})Y(b,z_{2})=
\varepsilon_{a,b} (z_{1}-z_{2})^{m}Y(b,z_{2})Y(a,z_{1}).
\end{eqnarray}

A vertex superalgebra $V$ is called a {\it vertex operator
superalgebra} (cf. [T], [FFR], [DL]) if
there is another distinguished vector $\omega$ of $V$ such that
\begin{eqnarray*}
&(V6)&[L(m),L(n)]=(m-n)L(m+n)+\frac{m^{3}-m}{12}\delta_{m+n,0}({\rm rank}V)\\
& &\mbox{for }m,n\in \Bbb{Z},\mbox{ where }
Y(\omega,z)=\sum_{n\in \Bbb{Z}}L(n)z^{-n-2},\;{\rm
rank}V\in \Bbb{C};\\
&(V7)&L(-1)=D,\; i.e., Y(L(-1)a,z)={d\over dz}Y(a,z)\;\;\;\mbox{ for any
}a\in V;\\
&(V8)&V\!=\!\oplus_{n\in {1\over 2}\Bbb{Z}}
V_{(n)}\mbox{ is }{1\over 2}\Bbb{Z}\mbox{-graded such that
 }L(0)|_{V_{(n)}}\!=\!nId_{V_{(n)}},\dim V_{(n)} \!<\!\infty,\;\;\;\;\\
& &\mbox{ and } V_{(n)}=0\;\;\mbox{ for }n\mbox{ sufficiently small. }
\end{eqnarray*}

Recall the following proposition from [FLM], [FHL], [DL] and [Li]:

\begin{propo}\label{DL, L1}
The Jacobi identity for
vertex superalgebra
can be equivalently replaced by the supercommutativity (2.23).
  \end{propo}

Let $V_{1}$ and $V_{2}$ be vertex superalgebras. A vertex superalgebra
 {\em homomorphism} from $V_{1}$ to $V_{2}$ is a linear map
$\sigma$ such that
\begin{eqnarray}
\sigma ({\bf 1})={\bf 1},\;\sigma (Y(a,z)b)=Y(\sigma
(a),z)\sigma (b)\;\;\;\mbox{for any } a,b\in V_{1}.
\end{eqnarray}
If both $V_{1}$ and $V_{2}$ are vertex operator superalgebras, a homomorphism
is required to map the Virasoro element of $V_{1}$ to the Virasoro element of
$V_{2}$.
An {\it automorphism} of a vertex (operator) superalgebra $V$ is
1-1 onto homomorphism from $V$ to $V$.
Let $\sigma$ be an automorphism of order $T$ of a  vertex
superalgebra $V$.
 Then $\sigma$ acts semisimply on $V$. Therefore
\begin{eqnarray}
V=V^{0}\oplus V^{1}\oplus \cdots \oplus V^{T-1}\end{eqnarray}
where $V^{k}$ is the eigenspace of $V$ for $\sigma$ with eigenvalue
${\rm exp}\left({2k\pi\sqrt{-1}\over T}\right)$. It is clear
that $V^{0}$ is a vertex subsuperalgebra of $V$ and all $V^{k}$
$(k=0,1,\cdots,T-1)$ are $V^{0}$-modules.

\begin{defi}\label{[D], [FFR], [FLM]}
  Let $(V,{\bf 1},D,Y)$ be a
vertex superalgebra with an automorphism $\sigma$
of order $T$. A $\sigma$-twisted $V$-module is a triple $(M,d,
Y_{M})$ consisting of a super vector
space $M$, a $\Bbb{Z}_{2}$-endomorphism $d$ of $M$ and a linear map
$Y_{M}(\cdot,z)$ from $V$ to
$({\rm End}M)[[z^{{1\over T}},z^{-{1\over T}}]]$ satisfying the
following conditions:

(M1)$\;\;\;\;$For any $a\in V,u\in M, a_{n}u=0$ for $n\in {1\over
T}\Bbb{Z}$ sufficiently large;

(M2)$\;\;\;\;Y_{M}({\bf 1},z)=Id_{M}$;

(M3)$\;\;\;\;\;[d,Y_{M}(a,z)]=Y_{M}(D(a),z)=\displaystyle{{d\over
dz}Y_{M}(a,z)}$ for any $a\in V$;

(M4)$\;\;\;\;$For any $\Bbb{Z}_{2}$-homogeneous $a,b\in V$, the following
$\sigma$-twisted
Jacobi identity holds:
\begin{eqnarray}
&
&\;\;\;z_{0}^{-1}\delta\!\left(\!\frac{z_{1}-z_{2}}{z_{0}}\!\right)\!Y_{M}(a,z_{1})
Y_{M}(b,z_{2})-
\varepsilon_{a,b}z_{0}^{-1}\delta\!\left(\!\frac{z_{2}-z_{1}}{-z_{0}}
\!\right)\!Y_{M}(b,z_{2})Y_{M}(a,z_{1})\;\;\nonumber\\
& &=z_{2}^{-1}\sum_{j=0}^{T-1}{1\over T}\delta\left(\left(\frac{z_{1}-z_{0}}
{z_{2}}\right)^{{1\over T}}\right)Y_{M}(Y(\sigma ^{j}a,z_{0})b,z_{2}).
\end{eqnarray}
\end{defi}

If $V$ is a vertex operator superalgebra, a $\sigma$-twisted
$V$-module for $V$ as a  vertex
superalgebra is called a  $\sigma$-twisted {\it weak} module
for $V$ as a vertex
operator superalgebra. A $\sigma$-twisted weak
$V$-module $M$ is said to be ${1\over 2T}\Bbb{Z}$-{\it graded}
if $M=\oplus_{n\in {1\over 2T}\Bbb{Z}}M(n)$ such that
$$(M5)\hspace{0.5cm}a_{n}M(r)\subseteq M(r+m-n-1)\;\;\mbox{
for }a\in V_{(m)}, m\in \Bbb{Z}, n\in {1\over T}\Bbb{Z}, r\in {1\over
2T}\Bbb{Z}.$$
A $\sigma$-twisted module $M$ for $V$ as a vertex superalgebra is called a
{\it $\sigma$-twisted  module} for $V$ as a vertex operator superalgebra  if
$M=\oplus_{\alpha \in \Bbb{C}}M_{(\alpha)}$ such that

$(M6)\;\;\;\;L(0)u=\alpha u\;\;\;\mbox{for }\alpha \in \Bbb{C},u\in
M_{(\alpha)};$

$(M7)\;\;\;\;\mbox{For any fixed }\alpha, M_{(\alpha +n)}=0\;\;\;\mbox{for
$n\in {1\over 2T}\Bbb{Z}$ sufficiently small;}$

$(M8)\;\;\;\;\dim M_{(\alpha )}<\infty\;\;\;\mbox{for any }\alpha\in
\Bbb{C}.$\\
It is clear that any $\sigma$-twisted module is a direct sum of ${1\over
2T}\Bbb{Z}$-graded $\sigma$-twisted modules.

\begin{rema}
 As in the untwisted case, the Virasoro algebra
relation follows from (M2)-(M4).
\end{rema}

If $a\in V^{k}$ is even or odd, the $\sigma$-twisted Jacobi identity
(2.26) becomes:
\begin{eqnarray}
& &z_{0}^{-1}\delta\!\left(\!\frac{z_{1}-z_{2}}{z_{0}}\!\right)\!Y_{M}(a,z_{1})
Y_{M}(b,z_{2})-\varepsilon_{a,b}z_{0}^{-1}\delta\!\left(\!
\frac{z_{2}-z_{1}}{-z_{0}}\!\right)\!
Y_{M}(b,z_{2})Y_{M}(a,z_{1})\nonumber\\
&=&z_{2}^{-1}\delta\!\left(\!\frac{z_{1}-z_{0}}
{z_{2}}\!\right)\left(\!\frac{z_{1}-z_{0}}
{z_{2}}\!\right)^{-{k\over T}}Y_{M}(Y(a,z_{0})b,z_{2}).
\end{eqnarray}
Setting $b={\bf 1}$, we get
\begin{eqnarray}
z_{2}^{-1}\delta\!\left(\!\frac{z_{1}-z_{0}}{z_{2}}\!\right)\!Y_{M}(a,z_{1})
=z_{2}^{-1}\delta\!\left(\!\frac{z_{1}-z_{0}}
{z_{2}}\right)\!\left(\!\frac{z_{1}-z_{0}}
{z_{2}}\!\right)^{-{k\over T}}\!Y_{M}(Y(a,z_{0}){\bf 1},z_{2}).
\end{eqnarray}
By taking ${\rm Res}_{z_{0}}z_{0}^{-1}$, we obtain
\begin{eqnarray}
z_{2}^{-1}\delta\left(\frac{z_{1}}{z_{2}}\right)Y_{M}(a,z_{1})
=z_{2}^{-1}\delta\left(\frac{z_{1}}
{z_{2}}\right)\left(\frac{z_{1}}
{z_{2}}\right)^{-{k\over T}}Y_{M}(a,z_{2}).
\end{eqnarray}
Therefore
\begin{eqnarray}
z^{{k\over T}}Y_{M}(a,z)\in ({\rm End}M)[[z,z^{-1}]]\;\;\;\mbox{ for any }
a\in V^{k}.
\end{eqnarray}
Taking ${\rm Res}_{z_{0}}$ of (2.27), we obtain the {\it twisted
supercommutator formula}:
\begin{eqnarray}
& &[Y_{M}(a,z_{1}),Y_{M}(b,z_{2})]_{\pm}\nonumber\\
&=:& Y_{M}(a,z_{1})Y_{M}(b,z_{2})-\varepsilon_{a,b}Y_{M}(b,z_{2})Y_{M}(a,z_{1})
\nonumber\\
&=&{\rm Res}_{z_{0}}z_{2}^{-1}\delta\left(\frac{z_{1}-z_{0}}
{z_{2}}\right)\left(\frac{z_{1}-z_{0}}
{z_{2}}\right)^{-{k\over T}}Y_{M}(Y(a,z_{0})b,z_{2})\nonumber\\
&=&\sum_{j=0}^{\infty}{1\over j!}\left(\left({\partial\over\partial
z_{2}}\right)^{j}z_{1}^{-1}\delta\left({z_{2}\over
z_{1}}\right)\left({z_{2}\over z_{1}}\right)^{{k\over
T}}\right)Y_{M}(a_{j}b,z_{2}).
\end{eqnarray}
Then we get the same
supercommutativity as the one in the untwisted case.
Multiplying (2.27) by $\displaystyle{\left(\frac{z_{1}-z_{0}}
{z_{2}}\right)^{{k\over T}}}$, then taking ${\rm Res}_{z_{1}}$, we
obtain the {\it twisted iterate formula}:
\begin{eqnarray}
Y_{M}(Y(a,z_{0})b,z_{2})={\rm Res}_{z_{1}}\left(\frac{z_{1}-z_{0}}
{z_{2}}\right)^{{k\over T}}\cdot X \end{eqnarray}
where
$$X\!=\!z_{0}^{-1}\delta\!\left(\!\frac{z_{1}-z_{2}}{z_{0}}\!\right)\!Y_{M}(a,z_{1})
Y_{M}(b,z_{2})-\varepsilon_{a,b}z_{0}^{-1}\delta\!\left(\!\frac{z_{2}-z_{1}}
{-z_{0}}\!\right)\!Y_{M}(b,z_{2})Y_{M}(a,z_{1}).$$
Let $r$ be a positive integer such that the
supercommutativity (2.23) holds. Then
\begin{eqnarray}
z_{0}^{r}X&=&(z_{1}-z_{2})^{r}X\nonumber \\
&=&z_{2}^{-1}\delta\!\left(\frac{z_{1}-z_{0}}
{z_{2}}\right)\left((z_{1}-z_{2})^{r}Y_{M}(a,z_{1})Y_{M}(b,z_{2})\right).
\end{eqnarray}
Thus
\begin{eqnarray}
& &\!\!z_{0}^{r}(z_{2}+z_{0})^{{k\over T}}Y_{M}(Y(a,z_{0})b,z_{2})\nonumber\\
&=&\!\!{\rm Res}_{z_{1}}z_{2}^{-1}\delta\!\left(\!\frac{z_{1}-z_{0}}
{z_{2}}\!\right)\!\left(\!\frac{z_{1}-z_{0}}{z_{2}}\!\right)^{{k\over
T}}\!(z_{2}+z_{0})^{{k\over T}}
\left((z_{1}-z_{2})^{r}Y_{M}(a,z_{1})Y_{M}(b,z_{2})\right)\nonumber \\
&=&\!\!{\rm Res}_{z_{1}}z_{1}^{-1}\delta\!\left(\!\frac{z_{2}+z_{0}}
{z_{1}}\!\right)\!\left(\!\frac{z_{2}+z_{0}}{z_{1}}\!\right)^{-{k\over T}}
\!(z_{2}+z_{0})^{{k\over T}}\!
\left((z_{1}-z_{2})^{r}Y_{M}(a,z_{1})Y_{M}(b,z_{2})\right)\nonumber \\
&=&\!\!{\rm Res}_{z_{1}}z_{1}^{-1}\delta\!\left(\frac{z_{2}+z_{0}}
{z_{1}}\right)
\left((z_{1}-z_{2})^{r}Y^{o}_{M}(a,z_{1})Y_{M}(b,z_{2})\right)\nonumber\\
&=&\!\!\left(z_{0}^{r}Y^{o}_{M}(a,z_{2}+z_{0})
Y_{M}(b,z_{2})\right),
\end{eqnarray}
where $Y^{o}_{M}(a,z)=z^{{k\over T}}Y_{M}(a,z)\in ({\rm End}M)[[z,z^{-1}]]$.
For any $u\in M$, let $k$ be a positive integer such that
$z^{k}Y^{o}_{M}(a,z)u\in M[[z]]$. Then we obtain
\begin{eqnarray}
& &(z_{2}+z_{0})^{{k\over
T}}(z_{0}+z_{2})^{k}z_{0}^{r}Y_{M}(Y(a,z_{0})b,z_{2})u\nonumber\\
&=&z_{0}^{r}(z_{0}+z_{2})^{k}Y^{o}_{M}(a,z_{0}+z_{2})
Y_{M}(b,z_{2})u.
\end{eqnarray}
Multiplying both sides by $z_{0}^{-r}$, we obtain the following twisted
associativity:
\begin{eqnarray}
& &(z_{0}+z_{2})^{k}(z_{2}+z_{0})^{{k\over
T}}Y_{M}(Y(a,z_{0})b,z_{2})u\nonumber \\
&=&(z_{0}+z_{2})^{k}Y^{o}_{M}(a,z_{0}+z_{2})Y_{M}(b,z_{2})u.
\end{eqnarray}

\begin{lemma}
 Let $\sigma$ be an automorphism  of order $T$
of a vertex superalgebra $V$. Then the $\sigma$-twisted Jacobi
identity (2.26) can
be equivalently replaced by the supercommutativity and the
$\sigma$-twisted  associativity (2.36).
\end{lemma}

{\bf Proof}. For any homogeneous $a,b\in V, u\in M$, let $m$ be a positive
integer such that both the supercommutativity and the twisted
associativity hold and that $z^{m}Y^{o}_{M}(a,z)u$ involves only
positive powers of $z$. Then
\begin{eqnarray*}
& &z_{0}^{m}z_{1}^{m}z_{0}^{-1}\delta\left(\frac{z_{1}-z_{2}}{z_{0}}
\right)Y_{M}(a,z_{1})Y_{M}(b,z_{2})u\nonumber\\
& &-z_{0}^{m}z_{1}^{m}\varepsilon_{a,b}
z_{0}^{-1}\delta\left(\frac{z_{2}-z_{1}}
{-z_{0}}\right)Y_{M}(b,z_{2})Y_{M}(a,z_{1})u\\
&=&z_{0}^{-1}\delta\left(\frac{z_{1}-z_{2}}{z_{0}}\right)\left(z_{1}^{m}
(z_{1}-z_{2})^{m}Y_{M}(a,z_{1})Y_{M}(b,z_{2})u\right)\\
& &-\varepsilon_{a,b}z_{0}^{-1}\delta\left(\frac{-z_{2}+z_{1}}{z_{0}}\right)
\left(z_{1}^{m}
(z_{1}-z_{2})^{m}Y_{M}(b,z_{2})Y_{M}(a,z_{1})u\right)\\
&=&z_{2}^{-1}\delta\left(\frac{z_{1}-z_{0}}{z_{2}}\right)\varepsilon_{a,b}
\left(z_{1}^{m}
(z_{1}-z_{2})^{m}Y_{M}(b,z_{2})Y_{M}(a,z_{1})u\right)\\
&=&z_{1}^{-1}\delta\left(\frac{z_{2}+z_{0}}{z_{1}}\right)\varepsilon_{a,b}
\left(z_{1}^{-{k\over T}}(z_{2}+z_{0})^{m}
z_{0}^{m}Y_{M}(b,z_{2})Y^{o}_{M}(a,z_{2}+z_{0})u\right)\\
&=&z_{1}^{-1}\delta\left(\frac{z_{2}+z_{0}}{z_{1}}\right)\varepsilon_{a,b}
\left(z_{1}^{-{k\over T}}z_{0}^{m}
(z_{0}+z_{2})^{m}Y_{M}(b,z_{2})Y^{o}_{M}(a,z_{0}+z_{2})u\right)\\
&=&z_{2}^{-1}\delta\left(\frac{z_{1}-z_{0}}{z_{2}}\right)
\left(z_{1}^{-{k\over T}}z_{0}^{m}
(z_{0}+z_{2})^{m}Y^{o}_{M}(a,z_{0}+z_{2})Y_{M}(b,z_{2})u\right)\\
&=&z_{2}^{-1}\delta\left(\frac{z_{1}-z_{0}}{z_{2}}\right)z_{1}^{-{k\over
T}}z_{0}^{m}(z_{0}+z_{2})^{m}(z_{2}+z_{0})^{{k\over T}}Y_{M}(Y(a,z_{0})b,
z_{2})u\\
&=&z_{0}^{m}z_{1}^{m}z_{2}^{-1}\delta\left(\frac{z_{1}-z_{0}}{z_{2}}\right)
\left(\frac{z_{2}+z_{0}}{z_{1}}\right)^{{k\over T}}Y_{M}(Y(a,z_{0})b,
z_{2})u. \end{eqnarray*}
Multiplying by $z_{0}^{-m}z_{1}^{-m}$, we obtain the $\sigma$-twisted
Jacobi identity. $\;\;\;\;\Box$

\begin{rema}
 Since the twisted associativity follows from
the supercommutativity and the
iterate formula (2.32), the supercommutativity and the
iterate formula (2.32) imply the $\sigma$-twisted Jacobi identity.
\end{rema}

Let $a,b\in V$ be
$\Bbb{Z}_{2}$-homogeneous such that $[Y_{M}(a,z_{1}),Y_{M}(b,z_{2})]_{\pm}=0$.
Then
\begin{eqnarray}
& &Y_{M}(Y(a,z_{0})b,z_{2})\nonumber\\
&=&{\rm Res}_{z_{1}}\left(\frac{z_{2}+z_{0}}{z_{1}}\right)^{-{k\over T}}
z_{1}^{-1}\delta\left(\frac{z_{2}+z_{0}}{z_{1}}\right)\left(Y_{M}(a,z_{1})
Y_{M}(b,z_{2})\right)\nonumber\\
&=&(z_{2}+z_{0})^{-{k\over T}}Y_{M}^{o}(a,z_{2}+z_{0})Y_{M}(b,z_{2})\nonumber\\
&=&Y_{M}(a,z_{2}+z_{0})Y_{M}(b,z_{2}).
\end{eqnarray}
Therefore
\begin{eqnarray}
Y_{M}(a_{-1}b,z_{2})=Y_{M}(a,z_{2})Y_{M}(b,z_{2}).
\end{eqnarray}
If $a\in V$ is $\Bbb{Z}_{2}$-homogeneous such that
$[Y_{M}(a,z_{1}),Y_{M}(a,z_{2})]_{\pm}=0$, then
\begin{eqnarray}
Y_{M}((a_{-1})^{N}{\bf 1},z)=Y_{M}(a,z)^{N}\;\;\;\mbox{for any }N\in
\Bbb{Z}_{+}.
\end{eqnarray}
This is Dong and Lepowsky's formula (13.70) in [DL] in the twisted case.
Then we have the following generalization of Proposition 13.16 in [DL].

\begin{propo}
 Let $V$ be a vertex superalgebra, let
$M$ be a $\sigma$-twisted $V$-module and let $a\in V$ be a
$\Bbb{Z}_{2}$-homogeneous
element such that $[Y_{M}(a,z_{1}),Y_{M}(a,z_{2})]_{\pm}=0$, so
that $Y(a,z)^{N}$ is well defined for $N\in \Bbb{Z}_{+}$. Then
$Y_{V}(a,z)^{N}=0$ implies $Y_{M}(a,z)^{N}=0$ on $M$. On the other
hand, the converse is true if $M$ is faithful.
\end{propo}

Let $\sigma$ be an automorphism of $V$ and let $M$ be a faithful
$\sigma$-twisted $V$-module. Let $a,b, u^{0},u^{1},\cdots,u^{n}\in V$.
It follows from Lemma 2.3 and the $\sigma$-twisted commutator formula
(2.31) that
\begin{eqnarray}
[Y_{M}(a,z_{1}),Y_{M}(b,z_{2})]_{\pm}
=\sum_{j=0}^{\infty}{1\over j!}\left(\left({\partial\over\partial
z_{2}}\right)^{j}z_{1}^{-1}\delta\!\left({z_{2}\over
z_{1}}\right)\!\left({z_{2}\over z_{1}}\right)^{{k\over
T}}\right)Y_{M}(u^{j},z_{2})
\end{eqnarray}
if and only if $a_{i}b=u^{j}$ for $0\le j\le n$, $a_{i}b=0$ for $j>n$.
Since $V$ is always a faithful $V$-module, we have:

\begin{lemma}
Let $V$ be a vertex superalgebra with an
automorphism of order $T$ and let $(M,Y_{M})$ be a faithful $\sigma$-twisted
$V$-module. Let $a,b,u^{0},\cdots,u^{n}$ be
$\Bbb{Z}_{T}$-homogeneous elements of $V$. Then
\begin{eqnarray}
[Y_{V}(a,z_{1}),Y_{V}(b,z_{2})]_{\pm}
=\sum_{j=0}^{\infty}{1\over j!}\left(\left({\partial\over\partial
z_{2}}\right)^{j}z_{1}^{-1}\delta\left({z_{2}\over
z_{1}}\right)\right)Y_{V}(u^{j},z_{2})
\end{eqnarray}
if and only if (2.40) holds.
\end{lemma}

\section{Local systems of twisted vertex operators}
Let us first define three basic categories $C,\;C^{o}$ and
$C_{\ell}$. The set of objects of the category $C$ is the
set of $\Bbb{Z}_{2}$-graded vector spaces and the set of
morphisms for any two objects is the set of all $\Bbb{Z}_{2}$-homomorphisms.
The set of objects of the category $C^{o}$
consists of all pairs $(M,d)$ where $M$ is an object of the category
$C$ and $d$ is an endomorphism of the object $M$. If
$(M_{i},d_{i})\;(i=1,2)$ are two objects of $C^{o}$, then a morphism is
a morphism $f$ from $M_{1}$ to $M_{2}$ such that $d_{2}f=fd_{1}$.
Finally,
$C_{\ell}$ is the category of $\Bbb{Z}_{2}$-graded vector spaces which are also
restricted $Vir$-modules with $\ell$ as its central charge and with the even
subspace and the odd subspace as submodules.

{\bf Throughout this section, $T$ will be a fixed positive integer.}

\begin{defi}
 Let $M$ be any
$\Bbb{Z}_{2}$-graded vector space.
A  $\Bbb{Z}_{T}$-{\it twisted weak vertex operator}
on $M$ is an element $a(z)=\sum_{n \in {1\over T}\Bbb{Z}}a_{n}z^{-n-1} \in
 ({\rm End}\:M)[[z^{{1\over T}},z^{-{1\over T}}]]$  such that
for any $u\in M, a_{n}u=0$ for $n\in {1\over T}\Bbb{Z}$ sufficiently large.
Let $(M,d)$ be an object of the category $C^{o}$.
A $\Bbb{Z}_{T}$-twisted weak vertex operator $a(z)$ on $M$ is
called a $\Bbb{Z}_{T}$-twisted weak vertex operator
on $(M,d)$ if it satisfies the following condition:
\begin{eqnarray}
[d,a(z)]=a'(z) \left(={d\over dz}a(z)\right).
\end{eqnarray}
Let $M$ be an object of $C_{\ell}$. A
$\Bbb{Z}_{T}$-twisted weak vertex
operator $a(z)$ on $(M,L(-1))$ is said to be of weight
$\lambda\in \Bbb{C}$ if it
satisfies the following condition:
\begin{eqnarray}
[L(0),a(z)]=\lambda a(z)+za'(z).
\end{eqnarray}
\end{defi}

Denote by $F(M,T)$ (resp. $F(M,d,T)$) the space of all
$\Bbb{Z}_{T}$-twisted weak vertex operators on $M$ (resp. $(M,d)$).
Set $\varepsilon=\exp \left(\frac{2\pi \sqrt{-1}}{T}\right)$.
Let $\sigma$ be the endomorphism of $({\rm End}M)[[z^{{1\over T}},z^{-{1\over
T}}]]$ defined by: $\sigma f(z^{{1\over T}})=f(\varepsilon^{-1} z^{{1\over
T}})$. It is clear that $\sigma$
restricted to $F(M,T)$ is an automorphism of order $T$. Then
\begin{eqnarray}
F(M,T)=F(M,T)^{0}\oplus F(M,T)^{1}\oplus\cdots\oplus F(M,T)^{T-1},
\end{eqnarray}
where $F(M,T)^{k}=\{f(z)\in F(M,T)|\sigma f(z)=\varepsilon^{k}f(z)\}$
for $0\le k\le T-1$. Therefore,
$F(M,T)$ is $\Bbb{Z}_{2}\times \Bbb{Z}_{T}$-graded. If $a(z)$ is a
$\Bbb{Z}_{T}$-twisted weak vertex operator on $M$ (resp. $(M,d)$),
then $\displaystyle{a'(z)={d\over dz}a(z)}$ is a
$\Bbb{Z}_{T}$-twisted weak vertex operator on $M$ (resp. $(M,d)$).
Then we have an endomorphism $\displaystyle{D={d\over dz}}$:
\begin{eqnarray}
D: F(M,T)\rightarrow {\rm F}(M,T);\;D\cdot a(z)={d\over dz}a(z)=a'(z).
\end{eqnarray}
It is clear that $D\sigma =\sigma D$.

\begin{defi}
 Two $\Bbb{Z}_{2}$-homogeneous
$\Bbb{Z}_{T}$-twisted weak vertex operators $a(z)$ and $b(z)$ are
said to be {\it mutually local} if there is a positive integer $n$ such that
\begin{eqnarray}
(z_{1}-z_{2})^{n}a(z_{1})b(z_{2})=(-1)^{|a(z)||b(z)|}(z_{1}-z_{2})^{n}
b(z_{2})a(z_{1}).
\end{eqnarray}
A $\Bbb{Z}_{T}$-twisted weak vertex operator $a(z)$ is called a
$\Bbb{Z}_{T}$-{\it twisted vertex
operator} if $a(z)$ is local with itself.
A $\Bbb{Z}_{2}$-graded subspace $A$ of $F(M,T)$ is said to be {\it
local} if any two $\Bbb{Z}_{2}$-homogeneous elements of $A$ are
mutually local. We define {\it a
local system} of $\Bbb{Z}_{T}$-twisted vertex operators on $M$
is a maximal $\Bbb{Z}_{2}$-graded local space of $\Bbb{Z}_{T}$-twisted
vertex operators on $M$.
\end{defi}

It is clear that any local system $A$ of
$F(M,T)$ has a natural
automorphism $\sigma$ such that $\sigma^{T}=Id_{A}$.

Replacing $n$ with $n+1$ in (3.5), then differentiating (3.5) with respect to
$z_{2}$ we have:

\begin{lemma}\label{Li3} [Li]
Let $M$ be a $\Bbb{Z}_{2}$-graded
vector space and let $a(z)$ and $b(z)$ be mutually local
$\Bbb{Z}_{T}$-twisted vertex operators on $M$.
Then $a(z)$ is local with $b'(z)$.
\end{lemma}

\begin{rema}
 Let $M$ be any $\Bbb{Z}_{2}$-graded vector
space. Then the identity operator $I(z)=Id_{M}$
is a $\Bbb{Z}_{T}$-twisted vertex operator which is mutually
local with any $\Bbb{Z}_{T}$-twisted weak vertex operator on $M$.
Therefore, any local system
of $\Bbb{Z}_{T}$-twisted vertex operators on $M$ contains
$I(z)$. It follows from Lemma 3.3 that any local system is closed
under the derivative operation $\displaystyle{D={d\over dz}}$.
\end{rema}

Using the same proof for the untwisted case in [Li] we get:

\begin{lemma}\label{Li2} [Li]
Let $M$ be a restricted $Vir$-module of
central charge
$\ell$. Then $L(z)=\sum_{n\in \Bbb{Z}}L(n)z^{-n-2}$ is a (local)
$\Bbb{Z}_{T}$-twisted vertex operator (of weight two) on $M$.
\end{lemma}

\begin{rema}
 Let $V$ be a vertex superalgebra with an
automorphism $\sigma$ of order $T$ and let $(M,Y_{M})$ be a $\sigma$-twisted
$V$-module. Then the image of $V$ under the linear map
$Y_{M}(\cdot,z)$ is a local subspace of $F(M,T)$.
\end{rema}

\begin{defi}
 Let $M$ be a $\Bbb{Z}_{2}$-graded vector
space and let $a(z)$ and $b(z)$ be mutually
local $\Bbb{Z}_{T}$-twisted vertex operators on $M$ such that
$a(z)\in F(M,T)^{k}$. Then for any integer $n$ we define
$a(z)_{n}b(z)$ as an element of $F(M,T)$ as follows:
\begin{eqnarray}
a(z)_{n}b(z)={\rm Res}_{z_{1}}{\rm
Res}_{z_{0}}\left(\frac{z_{1}-z_{0}}{z}\right)^{{k\over
T}}z_{0}^{n}\cdot X
\end{eqnarray}
where
$$X=z_{0}^{-1}\delta\!\left(\!\frac{z_{1}-z}{z_{0}}\!\right)\!a(z_{1})
b(z)-\varepsilon_{a,b}z_{0}^{-1}
\delta\!\left(\!\frac{z-z_{1}}{-z_{0}}\!\right)\!b(z)a(z_{1}),$$
or, equivalently, $a(z)_{n}b(z)$ is defined by:
\begin{eqnarray}
\sum_{n\in \Bbb{Z}}\left(a(z)_{n}b(z)\right)z_{0}^{-n-1}
={\rm Res}_{z_{1}}\left(\frac{z_{1}-z_{0}}{z}\right)^{{k\over T}}
\cdot X.
\end{eqnarray}
\end{defi}

\begin{rema}
 For any $u\in M$, we have:
\begin{eqnarray}
& &(a(z)_{n}b(z))u\nonumber\\
&=&\sum_{i=0}^{\infty}\left(\begin{array}{c}{k\over T}\\i\end{array}\right)
(-1)^{i}z_{1}^{{k\over T}-i}z^{{k\over T}}\cdot\nonumber\\
& &\cdot\left((z_{1}-z)^{n+i}a(z_{1})b(z)u-\varepsilon_{a,b}(-z+z_{1})^{n+i}
b(z)a(z_{1})u\right)\nonumber\\
&=&\sum_{i=0}^{N}\left(\begin{array}{c}{k\over T}\\i\end{array}\right)
(-1)^{i}z_{1}^{{k\over T}-i}z^{{k\over T}}\cdot\nonumber\\
& &\cdot\left((z_{1}-z)^{n+i}a(z_{1})b(z)u-\varepsilon_{a,b}(-z+z_{1})^{n+i}
b(z)a(z_{1})u\right).
\end{eqnarray}
Then it is clear that $(a(z)_{n}b(z))u\in M((z))$. Thus
$a(z)_{n}b(z)\in F(M,T)$.
It follows from the locality that $a(z)_{n} b(z)=0$ for
$n$ sufficiently large. Furthermore, for any nonnegative integer $n$, we have:
\begin{eqnarray}
& &z_{0}^{-1}\delta\left(\frac{z_{1}-z_{2}}{z_{0}}\right)\left(
\frac{z_{1}-z_{0}}
{z_{2}}\right)^{n}=z_{0}^{-1}\delta\left(\frac{z_{1}-z_{2}}{z_{0}}\right),\\
& &z_{0}^{-1}\delta\left(\frac{-z_{2}+z_{1}}{z_{0}}\right)
\left(\frac{z_{1}-z_{0}}{z_{2}}\right)^{n}=
z_{0}^{-1}\delta\left(\frac{-z_{2}+z_{1}}{z_{0}}\right).
\end{eqnarray}
Then
\begin{eqnarray}
\left(\frac{z_{1}-z_{0}}{z_{2}}\right)^{n}X=X.
\end{eqnarray}
Since $\displaystyle{\left(\frac{z_{1}-z_{0}}{z_{2}}\right)^{n}X}$
exists for any integer $n$, (3.11) is true for any integer $n$.
Then the right hand side of (3.7) depends only on $k\;mod\; T\Bbb{Z}$.
\end{rema}

The proof of the following proposition is similar to the proof in the
ordinary case [Li], which was given by Professor Chongying Dong.

\begin{propo}
 Let $a(z)$, $b(z)$ and $c(z)$ be
$\Bbb{Z}_{2}$-graded $\Bbb{Z}_{T}$-twisted weak vertex operators on $M$.
If both $a(z)$ and $b(z)$ are
local with $c(z)$, then
$a(z)_{n}b(z)$ is local with $c(z)$ for all $n\in \Bbb{Z}$.
\end{propo}

{\bf Proof}. Let $r$ be a positive integer such that $r+n>0$ and that the
following identities hold:
\begin{eqnarray*}
&
&(z_{1}-z_{2})^{r}a(z_{1})b(z_{2})=\varepsilon_{a,b}(z_{1}-z_{2})^{r}b(z_{2})a(z_{1}),\\
&
&(z_{1}-z_{2})^{r}a(z_{1})c(z_{2})=\varepsilon_{a,c}(z_{1}-z_{2})^{r}c(z_{2})a(z_{1}),\\
&
&(z_{1}-z_{2})^{r}b(z_{1})c(z_{2})=\varepsilon_{b,c}(z_{1}-z_{2})^{r}c(z_{2})b(z_{1}).
\end{eqnarray*}
Assume $a(z)\in F(M,T)^{i}$. Then we have
\begin{eqnarray}
a(z)_{n}b(z)
&=&{\rm
Res}_{z_{1}}\sum_{k=0}^{\infty}(-1)^{k}\left(\begin{array}{c}{i\over T}\\k
\end{array}\right)z_{1}^{{i\over T}-k}z_{2}^{-{i\over T}}A\nonumber\\
&=&{\rm
Res}_{z_{1}}\sum_{k=0}^{2r}(-1)^{k}\left(\begin{array}{c}{i\over T}\\k
\end{array}\right)z_{1}^{{i\over T}-k}z^{-{i\over T}}A,
\end{eqnarray}
where $A=(z_{1}-z)^{n+k}a(z_{1})b(z)
-\varepsilon_{a,b}(-z+z_{1})^{n+k}b(z)a(z_{1})$.
Since
\begin{eqnarray*}
& &(z-z_{3})^{4r}\left((z_{1}-z)^{n+k}a(z_{1})b(z)c(z_{3})
-\varepsilon_{a,b}(-z+z_{1})^{n+k}b(z)a(z_{1})c(z_{3})\right)\\
&=&\sum_{s=0}^{3r}\left(\begin{array}{c}3r\\s\end{array}\right)
(z-z_{1})^{3r-s}(z_{1}-z_{3})^{s}(z-z_{3})^{r}\cdot\\
& &\cdot \left((z_{1}-z)^{n+k}a(z_{1})b(z)c(z_{3})
-\varepsilon_{a,b}(-z+z_{1})^{n+k}b(z)a(z_{1})c(z_{3})\right)\\
&=&\sum_{s=r+1}^{3r}\left(\begin{array}{c}3r\\s\end{array}\right)
(z-z_{1})^{3r-s}(z_{1}-z_{3})^{s}(z-z_{3})^{r}\cdot\\
& &\cdot \left((z_{1}-z)^{n+k}a(z_{1})b(z)c(z_{3})
-\varepsilon_{a,b}(-z+z_{1})^{n+k}b(z)a(z_{1})c(z_{3})\right)\cdot\\
&=&\sum_{s=r+1}^{3r}\left(\begin{array}{c}3r\\s\end{array}\right)
(z-z_{1})^{3r-s}(z_{1}-z_{3})^{s}(z-z_{3})^{r}\cdot\\
& &\cdot \left((z_{1}-z)^{n+k}c(z_{3})a(z_{1})b(z)
-\varepsilon_{a,b}(-z+z_{1})^{n+k}c(z_{3})b(z)a(z_{1})\right)\\
&=&(z-z_{3})^{4r}\left((z_{1}-z)^{n+k}c(z_{3})a(z_{1})b(z)
-\varepsilon_{a,b}(-z+z_{1})^{n+k}c(z_{3})b(z)a(z_{1})\right),
\end{eqnarray*}
we have
\begin{eqnarray*}
(z-z_{3})^{4r}(a(z)_{n}b(z))c(z_{3})=(z-z_{3})^{4r}c(z_{3})(a(z)_{n}
b(z)).\;\;\;\;\Box
\end{eqnarray*}

Let $A$ be any local systems of $\Bbb{Z}_{T}$-twisted vertex
operators on $M$. Then it follows from Remark 3.8 and Proposition 3.9
that $a(z)_{n}b(z)\in A$ for $a(z)\in A^{k},b(z)\in A$. Then for any
$a(z)\in A^{k}, b(z)\in A$ we define
\begin{eqnarray}
Y_{A}(a(z),z_{0})b(z)=\sum_{n\in \Bbb{Z}}(a(z)_{n}b(z))z_{0}^{-n-1}.
\end{eqnarray}
Linearly extending the definition of $Y_{A}(\cdot,z_{0})$ to any
element of $A$, we obtain a linear map:
\begin{eqnarray}
Y_{A}(\cdot,z_{0}): A\rightarrow ({\rm End} A)[[z_{0},z_{0}^{-1}]];
a(z)\mapsto Y_{A}(a(z),z_{0}).
\end{eqnarray}
Also notice that $A$ contains the identity operator $I(z)$ and that $A$ is
closed under the derivative operation. For the rest of this section,
$A$ will be a fixed local system of $\Bbb{Z}_{T}$-twisted
vertex operators on $M$.

\begin{lemma}
For any $a(z)\in A$, we have
\begin{eqnarray}
Y_{A}(I(z),z_{0})a(z)=a(z),
\;Y_{A}(a(z),z_{0})I(z)=e^{z_{0}{\partial\over\partial
z}}a(z)=a(z+z_{0}).
\end{eqnarray}
\end{lemma}

{\bf Proof}. The proof of the first identity is the same as that in
the untwisted case. For the second one, without losing generality, we
may assume that $a(z)\in A^{k}$. Then
\begin{eqnarray}
& &\!\!Y_{A}(a(z),z_{0})I(z)\nonumber\\
&=&\!\!{\rm Res}_{z_{1}}\left(\!\frac{z_{1}-z_{0}}{z}\right)^{{k\over T}}
\!\left(z_{0}^{-1}\delta\!\left(\frac{z_{1}-z}{z_{0}}\right)a(z_{1})I(z)-
z_{0}^{-1}\delta\!\left(\frac{z-z_{1}}{-z_{0}}\right)I(z)a(z_{1})
\!\right)\nonumber\\
&=&{\rm Res}_{z_{1}}\left(\frac{z_{1}-z_{0}}{z}\right)^{{k\over T}}
z^{-1}\delta\left(\frac{z_{1}-z_{0}}{z}\right)a(z_{1})\nonumber\\
&=&{\rm Res}_{z_{1}}\left(\frac{z+z_{0}}{z_{1}}\right)^{-{k\over T}}
z_{1}^{-1}\delta\left(\frac{z+z_{0}}{z_{1}}\right) a(z_{1})\nonumber\\
&=&{\rm Res}_{z_{1}}(z+z_{0})^{-{k\over T}}z_{1}^{-1}\delta\left(\frac{z+z_{0}}
{z_{1}}\right)\left(z_{1}^{{k\over T}}a(z_{1})\right)\nonumber\\
&=&{\rm Res}_{z_{1}}(z+z_{0})^{-{k\over T}}z_{1}^{-1}\delta\left(\frac{z+z_{0}}
{z_{1}}\right)\left((z+z_{0})^{{k\over T}}a(z+z_{0})\right)\nonumber\\
&=&a(z+z_{0})\nonumber\\
&=&e^{z_{0}{\partial\over\partial z}}a(z).\;\;\;\;\Box
\end{eqnarray}

\begin{lemma}
For any $a(z),b(z)\in A$, we have:
\begin{eqnarray}
{\partial\over \partial z_{0}}Y_{A}(a(z),z_{0})b(z)=Y_{A}(a'(z),z_{0})b(z)=
Y_{A}(D(a(z)),z_{0})b(z).
\end{eqnarray}
\end{lemma}

{\bf Proof}. By linearity, we may assume that $a(z)\in F(M,T)^{k}$.
Then we have
\begin{eqnarray}
& &{\partial\over \partial z_{0}}Y_{A}(a(z),z_{0})b(z)\nonumber\\
&=&\!{\partial \over \partial z_{0}}{\rm
Res}_{z_{1}}\left(\frac{z_{1}-z_{0}}{z}\right)^{{k\over T}}\cdot\nonumber\\
& &\cdot \left( z_{0}^{-1}\delta\left(\frac{z_{1}-z}{z_{0}}\right)a(z_{1})
b(z)-\varepsilon_{a,b}z_{0}^{-1}
\delta\left(\frac{z-z_{1}}{-z_{0}}\right)b(z)a(z_{1})\right)\nonumber\\
&=&{\!\rm
Res}_{z_{1}}\left({\partial\over \partial z_{0}}
\left(\frac{z_{1}-z_{0}}{z}\right)^{{k\over T}}\right)\cdot\nonumber\\
& &\cdot \left( z_{0}^{-1}\delta\left(\frac{z_{1}-z}{z_{0}}\right)a(z_{1})
b(z)-\varepsilon_{a,b}z_{0}^{-1}
\delta\left(\frac{z-z_{1}}{-z_{0}}\right)b(z)a(z_{1})\right)\nonumber\\
& &+{\rm Res}_{z_{1}}\left(\frac{z_{1}-z_{0}}{z}\right)^{{k\over
T}}\cdot\nonumber\\
& &\cdot {\partial\over\partial z_{0}}\left(z_{0}^{-1}\delta \left(
\frac{z_{1}-z}{z_{0}}\right)a(z_{1})
b(z)+\varepsilon_{a,b}z_{0}^{-2}
\delta \left(\frac{z-z_{1}}{-z_{0}}\right)b(z)a(z_{1})\right)\nonumber\\
&=&\!-{\rm Res}_{z_{1}}\left({\partial\over \partial
z_{1}}\left({z_{1}-z_{0}\over z}\right)^{{k\over T}}\right)\cdot\nonumber\\
& &\cdot \left( z_{0}^{-1}\delta\left(\frac{z_{1}-z}{z_{0}}\right)a(z_{1})
b(z)-\varepsilon_{a,b}z_{0}^{-1}
\delta\left(\frac{z-z_{1}}{-z_{0}}\right)b(z)a(z_{1})\right)\nonumber\\
& &-{\rm Res}_{z_{1}}
\left({z_{1}-z_{0}\over z}\right)^{{k\over T}}\cdot\nonumber\\
& &\cdot \left(\left({\partial\over \partial z_{1}} z_{0}^{-1}\delta
\left(\frac{z_{1}-z}{z_{0}}\right)\right)a(z_{1})
b(z)-\varepsilon_{a,b}\left({\partial\over \partial z_{1}}z_{0}^{-1}
\delta\left(\frac{z-z_{1}}{-z_{0}}\right)\right)b(z)a(z_{1})\right)\nonumber\\
&=&\!{\rm Res}_{z_{1}}\!\left(\!{z_{1}-z_{0}\over z}\!\right)^{{k\over T}}\!
\left(\! z_{0}^{-1}\delta\!\left(\!\frac{z_{1}-z}{z_{0}}\!\right)\!a'(z_{1})
b(z)-\varepsilon_{a,b}z_{0}^{-1}
\delta\!\left(\!\frac{z-z_{1}}{-z_{0}}\!\right)\!b(z)a'(z_{1})\!\right)\nonumber\\
&=&\!Y_{A}(a'(z),z_{0})b(z). \;\;\;\;\Box
\end{eqnarray}

\begin{propo}
 Let $M$ be an object of category $C_{\ell}$ and let $A$ be a local system of
$\Bbb{Z}_{T}$-twisted vertex operators on $(M,L(-1))$.
Let $a(z)\in A^{j}, b(z)\in A^{k}$
of weights $\alpha$ and $\beta$, respectively. Then for any
integer $n$, $a(z)_{n}b(z)\in A^{j+k}$
has weight $(\alpha +\beta -n-1)$.
 \end{propo}

{\bf Proof.} It is enough to prove the following identities:
\begin{eqnarray}
& &\;\;\;[L(-1),Y_{A}(a(z),z_{0})b(z)]={\partial \over \partial
z}(Y_{A}(a(z),z_{0})b(z)),\\
& &\;\;\;[L(0),Y_{A}(a(z),z_{0})b(z)]\nonumber\\
& &=(\alpha +\beta)Y_{A}(a(z),z_{0})b(z)+z_{0}{\partial\over \partial
z_{0}}(Y_{A}(a(z),z_{0})b(z))+z{\partial\over \partial
z}(Y_{A}(a(z),z_{0})b(z)). \nonumber\\
& &\mbox{}
\end{eqnarray}
By definition we have:
\begin{eqnarray*}
& &{\partial\over \partial z}(Y_{A}(a(z),z_{0})b(z))\\
&=&{\partial\over \partial z}{\rm Res}_{z_{1}}
\left(\frac{z_{1}-z_{0}}{z}\right)^{{j\over T}}\cdot\\
& &\cdot\left(z_{0}^{-1}\delta\left(\frac{z_{1}-z}{z_{0}}\right)a(z_{1})b(z)
-\varepsilon_{a,b}z_{0}^{-1}\delta\left(\frac{z-z_{1}}{-z_{0}}\right)b(z)
a(z_{1})\right)\\
&=&{\rm Res}_{z_{1}}\left(-{j\over T}\right)z^{-1}\left(
\frac{z_{1}-z_{0}}{z}\right)^{{j \over T}}\cdot\\
& &\cdot\left(z_{0}^{-1}\delta\left(\frac{z_{1}-z}
{z_{0}}\right)a(z_{1})b(z)
-\varepsilon_{a,b}z_{0}^{-1}\delta\left(\frac{z-z_{1}}{-z_{0}}\right)b(z)
a(z_{1})\right)\\
& &+{\rm Res}_{z_{1}}\left(\frac{z_{1}-z_{0}}{z}\right)^{{j\over T}}\cdot\\
& &\cdot\left(-z_{0}^{-2}\delta '\left(\frac{z_{1}-z}
{z_{0}}\right)a(z_{1})b(z)
+\varepsilon_{a,b}z_{0}^{-2}\delta '\left(\frac{z-z_{1}}{-z_{0}}\right)b(z)
a(z_{1})\right)\\
& &+{\rm Res}_{z_{1}}\left(\frac{z_{1}-z_{0}}{z}\right)^{{j\over T}}\cdot\\
& &\cdot\left(\delta\left(\frac{z_{1}-z}{z_{0}}\right)a(z_{1})b'(z)
-\varepsilon_{a,b}z_{0}^{-1}\delta\left(\frac{z-z_{1}}{-z_{0}}\right)b'(z)
a(z_{1})\right)\\
&=&-{\rm Res}_{z_{1}}\left({\partial\over \partial z_{1}}
\left(\frac{z_{1}-z_{0}}{z}\right)^{{j\over T}}\right)\left(\frac{z_{1}-z_{0}}
{z}\right)\cdot\\
& &\cdot\left(z_{0}^{-1}\delta\left(\frac{z_{1}-z}
{z_{0}}\right)a(z_{1})b(z)
-\varepsilon_{a,b}z_{0}^{-1}\delta\left(\frac{z-z_{1}}{-z_{0}}\right)b(z)
a(z_{1})\right)\\
& &-{\rm Res}_{z_{1}}\left(\frac{z_{1}-z_{0}}{z}\right)^{{j\over T}}\cdot\\
& &\cdot\left(\left({\partial\over \partial z_{1}}z_{0}^{-1}\delta\left(\frac
{z_{1}-z}{z_{0}}\right)\right)a(z_{1})b(z)
-\left({\partial\over \partial z_{1}}z_{0}^{-1}\delta\left(\frac{z_{1}-z}
{z_{0}}\right)\right)b(z)a(z_{1})\right)\\
& &+{\rm Res}_{z_{1}}\left(\frac{z_{1}-z_{0}}{z}\right)^{{j\over T}}\cdot\\
& &\cdot\left(z_{0}^{-1}\delta\left(\frac{z_{1}-z}
{z_{0}}\right)a(z_{1})b'(z)
-\varepsilon_{a,b}z_{0}^{-1}\delta \left(\frac{z-z_{1}}{-z_{0}}\right)b'(z)
a(z_{1})\right)\\
&=&{\rm Res}_{z_{1}}\left(\frac{z_{1}-z_{0}}{z}\right)^{{j\over T}}\cdot\\
& &\cdot\left(z_{0}^{-1}\delta\left(\frac{z_{1}-z}
{z_{0}}\right)a'(z_{1})b(z)
-\varepsilon_{a,b}z_{0}^{-1}\delta \left(\frac{z-z_{1}}{-z_{0}}\right)b(z)
a'(z_{1})\right)\\
& &+{\rm Res}_{z_{1}}\left(\frac{z_{1}-z_{0}}{z}\right)^{{j\over T}}\cdot\\
& &\cdot\left(z_{0}^{-1}\delta\left(\frac{z_{1}-z}
{z_{0}}\right)a(z_{1})b'(z)
-\varepsilon_{a,b}z_{0}^{-1}\delta \left(\frac{z-z_{1}}{-z_{0}}\right)b'(z)
a(z_{1})\right)\\
&=&[L(-1),Y_{A}(a(z),z_{0})b(z)]
\end{eqnarray*}
and
\begin{eqnarray*}
& &[L(0),Y_{A}(a(z),z_{0})b(z)]\\
&=&{\rm Res}_{z_{1}}\left(\frac{z_{1}-z_{0}}{z}\right)^{{j\over T}}\cdot\\
& &\cdot\left(z_{0}^{-1}\delta\left(\frac{z_{1}-z_{2}}
{z_{0}}\right)[L(0),a(z_{1})b(z)]
-\varepsilon_{a,b}z_{0}^{-1}\delta \left(\frac{z-z_{1}}{-z_{0}}\right)
[L(0),b(z)a(z_{1})]\right)\\
&=&{\rm Res}_{z_{1}}\left(\frac{z_{1}-z_{0}}{z}\right)^{{j\over
T}}\cdot B
\end{eqnarray*}
where
\begin{eqnarray*}
B\!\!&=&\!\!\!z_{0}^{-1}\delta\!\left(\!\frac{z_{1}-z}
{z_{0}}\!\right)\!(a(z_{1})[L(0),b(z)]+[L(0),a(z_{1})]b(z))\\
& &\!\! -\varepsilon_{a,b}z_{0}^{-1}\delta\!\left(\frac{z_{1}-z}{z_{0}}\right)
(b(z)[L(0),a(z_{1})]+[L(0),b(z)]a(z_{1}))\\
&=&\!\!z_{0}^{-1}\delta\!\left(\!\frac{z_{1}-z}
{z_{0}}\!\right)\!\left(\beta a(z_{1})b(z)+za(z_{1})b'(z)+\alpha
a(z_{1})b(z)+z_{1}a'(z_{1})b(z)\right)\\
& &\!\!-z_{0}^{-1}\delta\!
\left(\!\frac{z_{1}-z}{z_{0}}\!\right)\!\left(\alpha b(z)a(z_{1})+z_{1}
b(z)a'(z_{1})+\beta b(z)a(z_{1})+zb'(z)a(z_{1})\right)\\
&=&\!\!z_{0}^{-1}\delta\!\left(\!\frac{z_{1}-z}
{z_{0}}\!\right)\!\left(\beta a(z_{1})b(z)+za(z_{1})b'(z)+\alpha
a(z_{1})b(z)+(z_{0}+z)a'(z_{1})b(z)\right)\\
& &\!\!-z_{0}^{-1}\delta\!
\left(\!\frac{z_{1}-z}{z_{0}}\!\right)\!\left(\alpha b(z)a(z_{1})+(z_{0}+z)
b(z)a'(z_{1})+\beta b(z)a(z_{1})+zb'(z)a(z_{1})\right)\!.
\end{eqnarray*}
Using Lemma 3.10 we obtain
\begin{eqnarray*}
& &[L(0),Y_{A}(a(z),z_{0})b(z)]\\
&=&(\alpha +\beta )Y_{A}(a(z),z_{0})b(z)+z{\partial\over \partial z}
(Y_{A}(a(z),z_{0})b(z))+z_{0}{\partial\over \partial z_{0}}(Y_{A}(a(z),z_{0})
b(z)).\\
& &\mbox{}\;\;\;\;\Box
\end{eqnarray*}

\begin{propo}
 Let $M$ be a $\Bbb{Z}_{2}$-graded
vector space and let $A$ be any local system of
$\Bbb{Z}_{T}$-twisted vertex operators on $M$. Then
$Y_{A}(a(z),z_{1})$ and $Y_{A}(b(z),z_{2})$ are mutually local
on $V$ for $a(z),b(z)\in A$.
\end{propo}

{\bf Proof}. By linearity, we may assume that $a(z)\in A^{j}$ and
$b(z)\in A^{k}$. Let $c(z)\in A$. Then we have:
\begin{eqnarray*}
& &Y_{A}(a(z),z_{3})Y_{A}(b(z),z_{0})c(z)\\
&=&{\rm Res}_{z_{1}}\left(\frac{z_{1}-z_{3}}{z}\right)^{{j\over T}}
z_{3}^{-1}\delta\left(\frac{z_{1}-z}{z_{3}}\right)
a(z_{1})(Y_{A}(b(z),z_{0})c(z))\\
& &-\varepsilon_{a,b}\varepsilon_{a,c}z_{3}^{-1}\delta\left(
\frac{-z+z_{1}}{z_{3}}\right)(Y_{A}(b(z),z_{0})c(z))a(z_{1})\\
&=&{\rm Res}_{z_{1}}{\rm Res}_{z_{4}}\left(\frac{z_{1}-z_{3}}{z}\right)^{
{j\over T}}\left(\frac{z_{4}-z_{0}}{z}\right)^{
{k\over T}}\cdot P
\end{eqnarray*}
where
\begin{eqnarray*}
P&=& z_{3}^{-1}\delta\left(\frac{z_{1}-z}{z_{3}}\right)z_{0}^{-1}
\delta\left(\frac{z_{4}-z}{z_{0}}\right)a(z_{1})b(z_{4})c(z)\\
& &- \varepsilon_{b,c}z_{3}^{-1}\delta\left(\frac{z_{1}-z}{z_{3}}\right)
z_{0}^{-1}
\delta\left(\frac{-z+z_{4}}{z_{0}}\right)a(z_{1})c(z)b(z_{4})\\
& &-\varepsilon_{a,b}\varepsilon_{a,c}
z_{3}^{-1}\delta\left(\frac{-z+z_{1}}{z_{3}}\right)z_{0}^{-1}
\delta\left(\frac{z_{4}-z}{z_{0}}\right)b(z_{4})c(z)a(z_{1})\\
& &+\varepsilon_{a,b}\varepsilon_{a,c}\varepsilon_{b,c}z_{3}^{-1}\delta
\left(\frac{-z+z_{1}}{z_{3}}\right)z_{0}^{-1}
\delta\left(\frac{-z+z_{4}}{z_{0}}\right)c(z)b(z_{4})a(z_{1}).
\end{eqnarray*}
Similarly, we have
\begin{eqnarray*}
& &Y_{A}(b(z),z_{0})Y_{A}(a(z),z_{3})c(z)\\
&=&{\rm Res}_{z_{1}}{\rm Res}_{z_{4}}\left(\frac{z_{1}-z_{3}}{z}\right)^{
{j\over T}}\left(\frac{z_{4}-z_{0}}{z}\right)^{
{k\over T}}\cdot Q
\end{eqnarray*}
where
\begin{eqnarray*}
Q&=& z_{3}^{-1}\delta\left(\frac{z_{1}-z}{z_{3}}\right)z_{0}^{-1}
\delta\left(\frac{z_{4}-z}{z_{0}}\right)b(z_{4})a(z_{1})c(z)\\
& &- \varepsilon_{b,c}z_{3}^{-1}\delta\left(\frac{-z+z_{1}}{z_{3}}\right)
z_{0}^{-1}
\delta\left(\frac{z_{4}-z}{z_{0}}\right)b(z_{4})c(z)a(z_{1})\\
& &-\varepsilon_{a,b}\varepsilon_{a,c}
z_{3}^{-1}\delta\left(\frac{z_{1}-z}{z_{3}}\right)z_{0}^{-1}
\delta\left(\frac{-z+z_{4}}{z_{0}}\right)a(z_{1})c(z)b(z_{4})\\
& &+\varepsilon_{a,b}\varepsilon_{a,c}\varepsilon_{b,c}z_{3}^{-1}\delta
\left(\frac{-z+z_{1}}{z_{3}}\right)z_{0}^{-1}
\delta\left(\frac{-z+z_{4}}{z_{0}}\right)c(z)a(z_{1})b(z_{4}).
\end{eqnarray*}
Let $k$ be any positive integer such that
\begin{eqnarray*}
(z_{1}-z_{4})^{k}a(z_{1})b(z_{4})=\varepsilon_{a,b}(z_{1}-z_{4})^{k}b(z_{4})a(z_{1}).
\end{eqnarray*}
Since
\begin{eqnarray}
& &(z_{3}-z_{0})^{k}z_{3}^{-1}\delta\left(\frac{z_{1}-z_{2}}{z_{3}}\right)
z_{0}^{-1}
\delta\left(\frac{z_{4}-z_{2}}{z_{0}}\right)\nonumber\\
&=&(z_{1}-z_{4})^{k}z_{3}^{-1}\delta\left(\frac{z_{1}-z_{2}}{z_{3}}\right)
z_{0}^{-1}\delta\left(\frac{z_{4}-z_{2}}{z_{0}}\right),
\end{eqnarray}
it is clear that
locality of $a(z)$ with $b(z)$ implies the locality of $Y_{A}(a(z),z_{1})$
with $Y_{A}(b(z),z_{2})$.$\;\;\;\;\Box$

\begin{theo}
 Let $M$ be a $\Bbb{Z}_{2}$-graded
vector space and let $A$ be any local system of
$\Bbb{Z}_{T}$-twisted vertex
operators on $M$. Then $V$ is a vertex superalgebra with an
automorphism $\sigma$ of order $T$ and $M$ is a
natural $\sigma$-twisted $A$-module in the sense
$Y_{M}(a(z),z_{1})=a(z_{1})$ for $a(z)\in A$.
\end{theo}

{\bf Proof.} It follows from  Proposition 3.9, Lemmas 3.10 and 3.11,
and Proposition 2.5 that $A$ is a vertex superalgebra. It is clear
that $\sigma$ is an automorphism of $A$. If we define
$Y_{M}(a(z),z_{0})=a(z_{0})$ for $a(z)\in A$, the
twisted vertex operator multiplication formula (3.6) becomes the
twisted iterate formula for $Y_{M}(\cdot,z_{0})$. Since the
supercommutativity (3.5) and iterate formula imply the twisted Jacobi
identity (Lemma 2.10), $Y_{M}(\cdot,z_{0})$ satisfies the twisted
Jacobi identity.
Therefore, $(M,Y_{M})$ is a $\sigma$-twisted $A$-module.$\;\;\;\;\Box$

Let $M$ be any $\Bbb{Z}_{2}$-graded vector space and let $S$ be
a set of $\Bbb{Z}_{2}$-homogeneous mutually local $\Bbb{Z}_{T}-$twisted
vertex operators on $M$. It follows from Zorn's lemma that there
exists a local system $A$ of $\Bbb{Z}_{T}$-twisted vertex operators
on $M$ containing $S$. Let $\langle S\rangle$ be the vertex
superalgebra generated by $S$ inside $A$ (which is a vertex
superalgebra). Since the ``multiplication'' (3.6) is canonical (which
is independent of the choice of the local system $A$), the vertex
superalgebra $\langle S\rangle$ is canonical, so that we may speak
about the vertex superalgebra generated by $S$ without
any confusion. Summarizing the previous argument, we have:

\begin{corol}
 Let $M$
be any $\Bbb{Z}_{2}$-graded vector space.
Then any set of mutually local $\Bbb{Z}_{T}$-twisted vertex operators
on $M$ generates a (canonical) vertex superalgebra with an
automorphism $\sigma$ of order $T$ such that $M$ is a $\sigma$-twisted
module for this vertex superalgebra.
\end{corol}

\begin{propo}\label{Li4} [Li]
 Let $M$ be an object of the category $C_{\ell}$ and let $V$ be a local system
of
$\Bbb{Z}_{T}$-twisted vertex operators on $(M,L(-1))$, which contains
$L(z)$. Then the vertex operator $L(z)$
is a Virasoro element of the vertex superalgebra $V$.
\end{propo}

{\bf Proof.} The proof is the same as the proof for untwisted case in
[Li].$\;\;\;\;\Box$

\begin{propo}
Let $V$ be a vertex superalgebra with an
automorphism $\sigma$ of order $T$. Then giving a $\sigma$-twisted
$V$-module $(M,d)$ is equivalent to giving a vertex superalgebra homomorphism
{}from $V$ to some local system of $\Bbb{Z}_{T}$-twisted vertex
operators on $(M,d)$.
\end{propo}

{\bf Proof.} Let $(M,d,Y_{M})$ be a $\sigma$-twisted $V$-module. Then
$\bar{V}=\{Y_{M}(a,z)|a\in V\}$ is a set of mutually local
$\Bbb{Z}_{T}$-twisted vertex operators on $(M,d)$. By Zorn's lemma, there exist
local systems which contain $\bar{V}$. Let $A$ be any one of those.
Then by definition $Y_{M}(\cdot,z)$ is a linear map from $V$ to $A$.
For homogeneous $a\in V^{j},b\in V$ we have:
\begin{eqnarray}
& &\!Y_{M}(\cdot,z)(Y_{V}(a,z_{0})b)\nonumber\\
&=&\!Y_{M}(Y_{V}(a,z_{0})b,z)\nonumber\\
&=&\!{\rm Res}_{z_{1}}\left(\!\frac{z_{1}-z_{0}}
{z_{2}}\!\right)^{{j\over T}}\cdot \nonumber\\
&
&\!\left(\!z_{0}^{-1}\delta\!\left(\!\frac{z_{1}-z}{z_{0}}\!\right)\!Y_{M}(a,z_{1})
Y_{M}(b,z)-\varepsilon_{a,b}z_{0}^{-1}\delta\!\left(\!\frac{z-z_{1}}
{-z_{0}}\!\right)\!Y_{M}(b,z)Y_{M}(a,z_{1})\!\right)\nonumber\\
&=&\!Y_{A}(Y_{M}(a,z),z_{0})Y_{M}(b,z)\nonumber\\
&=&\!Y_{A}(Y_{M}(\cdot,z)(a),z_{0})Y_{M}(\cdot,z)(b).
\end{eqnarray}
By definition we have:
\begin{eqnarray}
Y_{M}(\cdot,z)({\bf 1}_{V})=Y_{M}({\bf 1},z)=Id_{M}=I(z).
\end{eqnarray}
Thus $Y_{M}(\cdot,z)$ is a vertex superalgebra homomorphism from $V$
to $A$. Conversely, let $A$ be a local system of
$\Bbb{Z}_{T}$-twisted vertex operators on $M$ and let $Y_{M}(\cdot,z)$ be a
vertex superalgebra homomorphism from $V$ to $A$. Then
\begin{eqnarray}
Y_{M}({\bf 1},z)=Y_{M}(\cdot,z)({\bf 1})=I(z)=Id_{M}.
\end{eqnarray}
Since $A$ is a local system, for any $a,b\in V$, $Y_{M}(a,z)$ and
$Y_{M}(b,z)$ satisfy the supercommutativity (2.23). By reversing (3.22) we
obtain the iterate formula for $Y_{M}(\cdot,z)$. Since the twisted
Jacobi identity follows from the supercommutativity and the iterate
formula, $Y_{M}(\cdot,z)$ satisfies the twisted Jacobi identity.
Therefore, $(M,Y_{M})$ is a $\sigma$-twisted $V$-module.
$\;\;\;\;\Box$

\section{ Twisted modules for vertex operator superalgebras
associated to some infinite-dimensional Lie superalgebras}
In this section, we shall use the machinery we built in Section 3 to
study  twisted modules for vertex operator (super)algebras associated
to the Neveu-Schwarz algebra and  an affine
Lie superalgebra.

\subsection{Twisted modules for the Neveu-Schwarz vertex operator
superalgebra}
Let us first recall the Neveu-Schwarz algebra
(cf. [FFR], [KW], [T]). The Neveu-Schwarz algebra is the Lie superalgebra
\begin{eqnarray}
NS=\oplus_{m\in \Bbb{Z}}{\bf C}L(m)\oplus \oplus_{n\in \Bbb{Z}}
\Bbb{C}G(n+\frac{1}{2})\oplus \Bbb{C}c
\end{eqnarray}
with the following commutation relations:
\begin{eqnarray}
& &[L(m),L(n)]=(m-n)L(m+n)+\frac{m^{3}-m}{12}\delta_{m+n,0}c,\\
& &[L(m),G(n+\frac{1}{2})]=\left(\frac{m}{2}-n-\frac{1}{2}\right)G(m+n+
\frac{1}{2}),\\
& &[G(m+\frac{1}{2}),G(n-\frac{1}{2})]_{+}=2L(m+n)+\frac{1}{3}m(m+1)
\delta_{m+n,0}c,\\
& &[L(m),c]=0,\;[G(n+\frac{1}{2}),c]=0.
\end{eqnarray}
Or, equivalently
\begin{eqnarray}
&
&\!\![L(z_{1}),L(z_{2})]\!=\!z_{1}^{-1}\delta\!\left(\!\frac{z_{2}}{z_{1}}\!\right)\!
L'(z_{2})+2z_{1}^{-2}\delta'\!\left(\!\frac{z_{2}}{z_{1}}\!\right)\!L(z_{2})
\!+\!\frac{\ell}{12}z_{1}^{-4}\delta^{(3)}\!\left(\!\frac{z_{2}}{z_{1}}\!\right)\!,
\;\;\;\;\;\;\\
& &\!\![L(z_{1}),G(z_{2})]
=z_{1}^{-1}\delta\!\left(\frac{z_{2}}{z_{1}}\right){\partial
\over\partial z_{2}}G(z_{2})+\frac{3}{2}\left({\partial\over\partial z_{2}}
z_{1}^{-1}\delta\!\left(\frac{z_{2}}{z_{1}}\right)\right)G(z_{2}),
\end{eqnarray}
and
\begin{eqnarray}
[G(z_{1}),G(z_{2})]_{+}
=2z_{1}^{-1}\delta\left(\frac{z_{2}}{z_{1}}\right)L(z_{2})+\frac{1}{3}c
\left({\partial\over\partial z_{2}}\right)^{2}z_{1}^{-1}\delta\left(\frac
{z_{2}}{z_{1}}\right)
\end{eqnarray}
where
$L(z)=\sum_{m\in \Bbb{Z}}L(m)z^{-m-2}$ and
$G(z)=\sum_{n\in \Bbb{Z}}G(n+\frac{1}{2}) z^{-n-2}$.

By definition, $NS^{(0)}=\oplus_{m\in \Bbb{Z}}\Bbb{C}L(m)\oplus \Bbb{C}c$ and
$NS^{(1)}= \oplus_{n\in \Bbb{Z}}\Bbb{C}G(n+\frac{1}{2})$.
$NS=\oplus_{n\in {1\over 2}\Bbb{Z}}NS_{n}$ is a ${1\over
2}\Bbb{Z}$-graded Lie superalgebra by defining
\begin{eqnarray}
\deg L(m)=m,\;\deg c=0,\;\deg G(n)=n\;\;\;\mbox{for }m\in \Bbb{Z},n\in
{1\over 2}+\Bbb{Z}.
\end{eqnarray}
We have a triangular decomposition $NS=NS_{+}\oplus
NS_{0}\oplus NS_{+}$ where
\begin{eqnarray}
NS_{\pm}=\sum_{n=1}^{\infty}(\Bbb{C}L(\pm n)+\Bbb{C}G(\pm n\mp{1\over 2})),\;
NS_{0}=\Bbb{C}L(0)\oplus \Bbb{C}c.
\end{eqnarray}

Let $\bar{M}(c,0)=M(c,0)/\langle G(-{1\over 2}){\bf 1}\rangle$, where
$M(c,0)$ is the Verma module over $NS$ of central charge $c$ and
$\langle G(-{1\over 2}){\bf 1}\rangle$ is the submodule generated by
$G(-{1\over 2}){\bf 1}$. It has been proved ([KW], [Li]) that
$\bar{M}(c,0)$ has a natural vertex
operator superalgebra structure and any restricted $NS$-module of
charge $c$ is a weak $\bar{M}(c,0)$-module.

Let $V$ be any vertex superalgebra. Then the linear map $\sigma$:
$V\rightarrow V;\sigma (a+b)=a-b$ for $a\in V^{(0)}, b\in V^{(1)}$, is
an automorphism of $V$, which was called the {\it canonical
automorphism } [FFR]. It is
easy to see that $Aut\bar{M}(c,0)=\Bbb{Z}_{2}=\langle
\sigma\rangle$.

The Ramond algebra ${\bf R}$ [GSW] is the Lie superalgebra with a basis $L(m),
F(n),c$ for $m,n\in \Bbb{Z}$ and with the following defining relations:
\begin{eqnarray}
& &[L(m),L(n)]=(m-n)L(m+n)+\frac{m^{3}-m}{12}\delta_{m+n,0}c,\\
& &[L(m), F(n)]=({m\over 2}-n)F(m+n),\\
& &[F(m),F(n)]_{+}=2L(m+n)+{1\over 3}(m^{2}-{1\over 4})\delta_{m+n,0}c,\\
& &[L(m),c]=[F(m),c]=0
\end{eqnarray}
for $m,n\in \Bbb{Z}$. Set
\begin{eqnarray}
L(z)=\sum_{m\in \Bbb{Z}}L(m)z^{-m-2},\;
\; F(z)=\sum_{n\in \Bbb{Z}}F(n)z^{-n-{3\over 2}}.
\end{eqnarray}
Then the defining relations (4.12) and (4.13) are equivalent to the
following identities:
\begin{eqnarray}
& &\!\![L(z_{1}), F(z_{2})]\!=\!z_{1}^{-1}\delta\!\left(\!{z_{2}\over
z_{1}}\!\right)\!F(z_{2})
\!+\!{3\over 2}\!\left(\!{\partial\over\partial z_{2}}z_{1}^{-1}\delta\!\left(
\!{z_{2}\over z_{1}}\!\right)\!\left(\!{z_{2}\over z_{1}}\!\right)^{{1\over
2}}\right)\!
F(z_{2}),\;\;\;\\
& &\!\![F(z_{1}),F(z_{2})]_{+}\!=\!2\delta\!\left(\!{z_{2}\over
z_{1}}\!\right)\!\left(\!
{z_{2}\over z_{1}}\!\right)^{{1\over 2}}L(z_{1})+{2\over
3}\left(\!{\partial\over\partial
z_{2}}\!\right)^{2}\!\left(\!\delta\!\left(\!{z_{2}\over
z_{1}}\!\right)\!\left(\!{z_{2}\over z_{1}}\!\right)^{{1\over 2}}\right)\!.
\end{eqnarray}

Let $M$ be any restricted module for the Ramond algebra or the
twisted Neveu-Schwarz algebra $\Bbb{R}$
with central charge $c$. Then $\{ L(z), F(z), I(z)\}$ is a set of
mutually local homogeneous $\Bbb{Z}_{2}$-twisted vertex
operators on $M$.
Let $V$ be the vertex superalgebra generated by $\{L(z),F(z),I(z)\}$.
Then $M$ is a faithful $\sigma$-twisted $V$-module.
By Lemma 2.12, $Y_{V}(L(z),z_{1})$ and  $Y_{V}(F(z),z_{1})$ satisfy
the Neveu-Schwarz relations (4.6) and (4.7).
Then $V$ is a lowest weight module for $NS$ under
the linear map: $L(m)\mapsto L(z)_{m+1}, G(n+{1\over 2}) \mapsto
F(z)_{n}, c\mapsto c$. Then $V$ is a vertex operator superalgebra,
which is a quotient algebra of $\bar{M}(c,0)$. Consequently, $M$ is a
weak $\sigma$-twisted module for $\bar{M}(c,0)$. Therefore, we have proved:

\begin{propo}
 Let $c$ be any complex number. Then any
restricted module for the Ramond algebra $R$ of central charge $c$ is
a weak $\sigma$-twisted $\bar{M}(c,0)$-module.
 \end{propo}

\subsection{Twisted modules for vertex operator superalgebras
associated to affine Lie superalgebras}
Let $({\bf g},B)$ be a pair consisting of a finite-dimensional
Lie superalgebra ${\bf g}={\bf g}^{(0)}\oplus {\bf g}^{(1)}$ such that
$[{\bf g}^{(1)},{\bf g}^{(1)}]_{+}=0$ and a
nondegenerate symmetric bilinear form $B$ such that
\begin{eqnarray}
B({\bf g}^{(0)},{\bf g}^{(1)})=0,\; B([a,u],v)=-B(u,[a,v])\;\;\;\mbox{for }a\in
{\bf g}^{(0)}, u,v\in {\bf g}.
\end{eqnarray}
This amounts to having a finite-dimensional Lie algebra ${\bf g}^{(0)}$
with a nondegenerate symmetric invariant bilinear form
$B_{0}(\cdot,\cdot)$ and a finite-dimensional ${\bf g}^{(0)}$-module
${\bf g}^{(1)}$ with a nondegenerate symmetric bilinear form
$B_{1}(\cdot,\cdot)$ such that
\begin{eqnarray}
B_{1}(au,v)=-B_{1}(u,av)\;\;\;\mbox{ for any }a\in {\bf g}^{(0)}, u,v\in
{\bf g}^{(1)}.
\end{eqnarray}
The affine Lie superalgebra $\tilde{{\bf g}}$ is defined
to be
$\Bbb{C}[t,t^{-1}]\otimes {\bf g}\oplus \Bbb{C}c$ with the following
defining relations:
\begin{eqnarray}
& &[a_{m},b_{n}]=[a,b]_{m+n}+m\delta_{m+n,0}B(a,b)c,\\
& &[a_{m},u_{n}]=-[u_{n},a_{m}]=(au)_{m+n},\\
& &[u_{m},v_{n}]_{+}=\delta_{m+n+1,0}B(u,v) c,\\
& &[c, \tilde{{\bf g}}]=0
\end{eqnarray}
for $a,b\in {\bf g}^{(0)}, u,v\in {\bf g}^{(1)}$ and
$m,n\in \Bbb{Z}$, where $x_{m}$ stands for $t^{m}\otimes x$.
Or, equivalently
\begin{eqnarray}
& &[a(z_{1}),b(z_{2})]=z_{1}^{-1}\delta\left(\frac{z_{2}}{z_{1}}\right)
[a,b](z_{2})+
z_{1}^{-2}\delta'\left(\frac{z_{2}}{z_{1}}\right)
B(a,b)c,\\
& &[a(z_{1}),u(z_{2})]=z_{1}^{-1}\delta\left(\frac{z_{2}}{z_{1}}\right)
(au)(z_{2}),\\
& &[u(z_{1}),v(z_{2})]_{+}=B(u,v)
z_{1}^{-1}\delta\left(\frac{z_{2}}{z_{1}}\right)c,\\
& &[x(z),c]=0
\end{eqnarray}
for any $a,b\in {\bf g}^{(0)}, u,v\in {\bf g}^{(1)}, x\in {\bf g}, m,n\in
\Bbb{Z}$, where
$x(z)=\sum_{n\in \Bbb{Z}}x_{n}z^{-n-1}$ for $x\in {\bf g}$.
The even subspace and the odd
subspace are:
\begin{eqnarray}
\tilde{{\bf g}}^{(0)}=\Bbb{C}[t,t^{-1}]\otimes {\bf g}^{(0)}\oplus \Bbb{C},\;\;
\tilde{{\bf g}}^{(1)}={\bf C}[t,t^{-1}]\otimes {\bf g}^{(1)}.
\end{eqnarray}
Furthermore, $\tilde{{\bf g}}$ is a ${1\over 2}\Bbb{Z}$-graded Lie superalgebra
with the degree defined as follows:
\begin{eqnarray}
\deg a_{m}=-m, \;\deg u_{n}=-n+{1\over 2},\; \deg c=0
\end{eqnarray}
for any $a\in {\bf g}^{(0)}, u\in {\bf g}^{(1)}, m,n\in \Bbb{Z}$.
We have a triangular decomposition $\tilde{{\bf g}}=N_{+}\oplus
N_{0}\oplus N_{-}$, where
\begin{eqnarray}
N_{-}=t\Bbb{C}[t]\otimes {\bf g},\; N_{+}=t^{-1}\Bbb{C}[t^{-1}]\otimes
{\bf g}^{(0)}\oplus \Bbb{C}[t^{-1}]\otimes
{\bf g}^{(1)},\;N_{0}={\bf g}^{(0)}\oplus \Bbb{C}c.
\end{eqnarray}
Let $P=N_{-}+N_{0}$ be a parabolic subalgebra.
For any ${\bf g}^{(0)}$-module $U$ and any complex
number $\ell$, denote by $M_{({\bf g},B)}(\ell,U)$ the generalized Verma
module or Weyl
module for $\tilde{{\bf g}}$ with $c$ acting as scalar $\ell$. Namely,
$M_{({\bf g},B)}(\ell,U)=U(\tilde{{\bf g}})\otimes_{U(P)}U$.

It has been proved ([KW],[Li]) that for any complex number $\ell$,
$M_{({\bf g},B)}(\ell,\Bbb{C})$ has a natural vertex superalgebra
structure and that any restricted $\tilde{{\bf g}}$-module $M$ of level
$\ell$ is a module for the vertex superalgebra
$M_{({\bf g},B)}(\ell,\Bbb{C})$. Furthermore, if ${\bf g}^{(0)}$ is simple and
$\ell$ is not zero, then $M_{({\bf g},B)}(\ell,\Bbb{C})$ is a vertex
operator superalgebra. For any $a\in {\bf g}$, we may identify $a$
with $a_{-1}{\bf 1}\in M_{({\bf g},B)}(\ell,\Bbb{C})$, so that we may
consider ${\bf g}$ as a subspace of $M_{({\bf g},B)}(\ell,\Bbb{C})$.

{\it An automorphism of $({\bf g},B)$} is an automorphism $\sigma$ of the Lie
superalgebra ${\bf g}$, which preserves the even and the odd subspaces
and the bilinear form $B$. Let $\sigma$ be an automorphism of
$({\bf g},B)$ of order $T$. Then we have:
\begin{eqnarray}
{\bf g}={\bf g}_{0}\oplus {\bf g}_{1}\oplus\cdots\oplus {\bf g}_{T-1}\;\;\;
\mbox{where }{\bf g}_{j}=\{u\in {\bf g}|\sigma u=\varepsilon^{j}u\}.
\end{eqnarray}
The twisted Lie superalgebra $\tilde{{\bf g}}[\sigma]$ is defined to be
$\displaystyle{\tilde{{\bf g}}[\sigma]=\oplus_{j=0}^{T-1}{\bf
g}_{j}\otimes t^{{j\over T}}\Bbb{C}[t,t^{-1}]\oplus \Bbb{C}c}$ with
the following defining relations:
\begin{eqnarray}
& &\!\![t^{m+{j\over T}}\!\otimes\! a,t^{n+{k\over T}}\!\otimes\! b]
\!=\!\!t^{m+n+\frac{j+k}{T}}\!\otimes\! [a,b]\!+\!B(a,b)\!\left(\!m+{j\over
T}\!\right)\!
\delta_{m+n+\frac{j+k}{T},0}c,\;\;\;\;\;\\
& &\!\![t^{m+{j\over T}}\otimes u,t^{n+{k\over T}}\otimes v]
=B(a,b)\delta_{m+n+1+\frac{j+k}{T},0}c,\\
& &\!\![c,\tilde{{\bf g}}]=0
\end{eqnarray}
for $a\in ({\bf g}^{(0)})_{j}, b\in {\bf g}_{k}, u\in ({\bf g}^{(1)})_{j},
v\in ({\bf g}^{(1)})_{k}$. For $a\in {\bf g}_{j}$, set
\begin{eqnarray}
a(z)=\sum_{n\in \Bbb{Z}}(t^{n+{j\over T}}\otimes a)z^{-n-1-{j\over
T}}\in ({\rm End}M)[[z^{{1\over T}},z^{-{1\over T}}]].
\end{eqnarray}
Then for $a\in ({\bf g}^{(0)})_{j},b\in {\bf g}_{k}, u\in ({\bf
g}^{(1)})_{j}, v\in {\bf g}^{(1)}$, we have:
\begin{eqnarray}
& &\![a(z_{1}),b(z_{2})]\nonumber \\
&=&\!\sum_{m,n\in \Bbb{Z}}t^{m+n+{j\over T}+{k\over T}}\otimes
[a,b]z_{1}^{-m-1-{j\over T}}z_{2}^{-n-1-{k\over
T}}\nonumber\\
& &+\!c(m+{j\over T})B(a,b)
\delta_{m,-n}\delta_{j,-k}z_{1}^{-m-1-{j\over
T}}z_{2}^{-n-1-{k\over T}}\nonumber \\
&=&\!z_{1}^{-1}\left({z_{2}\over z_{1}}\right)^{{j\over
T}}\delta\!\left({z_{2}\over z_{1}} \right)[a,b](z_{2})
+B(a,b) {\partial\over\partial z_{2}}\!\left(z_{1}^{-1}
\delta\!\left({z_{2}\over z_{1}}\right)\!
\left({z_{2}\over z_{1}}\right)^{{j\over T}}\right)c,\;\;\;\\
& &[u(z_{1}),v(z_{2})]_{+}=B(u,v)
z_{1}^{-1}\delta\left({z_{2}\over z_{1}}\right)
\left({z_{2}\over z_{1}}\right)^{{j\over T}}c.
\end{eqnarray}

For any $\tilde{{\bf g}}$-module $M$, we may consider $a(z)$  as an
element of $(\!{\rm End}M\!)[[z^{{1\over T}},z^{-{1\over T}}]]$ for $a\in
{\bf g}$.
A $\tilde{{\bf g}}[\sigma]$-module $M$ is said to be {\it restricted} if
for any $u\in M$, $\left(t^{n+{j\over T}}\otimes {\bf
g}_{j}\right)\!u\!=0$ for $n$ sufficiently large. Then a $\tilde{{\bf
g}}[\sigma]$-module $M$ is restricted if and only if every $a(z)$ for
$a\in {\bf g}$ is a weak  $\Bbb{Z}_{T}$-twisted vertex
operators on $M$. Let $M$ be a restricted $\tilde{{\bf g}}[\sigma]$-module of
level $\ell$. It follows from (4.36), (4.37) and Lemma 2.2 that
$\{a(z)| a\in {\bf g}\}$ is a local subspace of $F(M,T)$. Let
$V$ be the vertex superalgebra generated by
$\{a(z)|a\in {\bf g}\}$ from the identity operator $I(z)$. Then $M$ is a
faithful $\sigma$-twisted $V$ module. It follows from (4.36), (4.37)
and Lemma 2.13 that $V$ is a $\tilde{{\bf g}}$-module where $a_{m}$
for $a\in {\bf g},m\in \Bbb{Z}$ is represented by $a(z)_{m}$. It is
easy to see that $V$ is a lowest weight $\tilde{{\bf g}}$-module with
a lowest weight vector $I(z)$, so that $V$ is a quotient vertex
superalgebra of $M_{({\bf g},B)}(\ell,\Bbb{C})$. Through the natural
vertex superalgebra homomorphism, $M$ becomes a $\sigma$-twisted
$M_{({\bf g},B)}(\ell,\Bbb{C})$-module. Thus we have proved:

\begin{theo}
 Let $({\bf g},B)$ be a pair given as before,
let $\sigma$ be any automorphism of order
$T$ of $M_{({\bf g},B)}(\ell,\Bbb{C})$
or $({\bf g},B)$ and let $\ell$ be any complex number. Then
any restricted $\tilde{{\bf g}}[\sigma]$-module of level $\ell$ is a
$\sigma$-twisted $M_{({\bf g},B)}(\ell,\Bbb{C})$-module.
\end{theo}

Next we shall consider two special cases.

{\bf Case 1.} Let ${\bf g}^{(0)}=0$ and let ${\bf g}^{(1)}=A$ be an
$n$-dimensional vector space with a nondegenerate symmetric bilinear
form $B$. It has been proved ([KW], [Li]) that if $\ell\ne 0$,
$M_{({\bf g}_{1},B)}(\ell,{\bf
C})$ ($=\Bbb{F}^{n}$ [KW], [Li], or $=CM(\Bbb{Z}+{1\over 2})$ [DM1],
[FFR]) is a rational vertex operator superalgebra which
has only one irreducible module up to equivalence.

Let $\sigma$ be an automorphism of $(A,B)$ of order $T$. Then
\begin{eqnarray}
A=A_{0}\oplus A_{2}\oplus\cdots\oplus A_{T-1}
\end{eqnarray}
where $A_{j}=\{u\in A|\sigma (u)=\varepsilon^{j}u\}$. Set
\begin{eqnarray}
\tilde{A}[\sigma]
=\oplus_{j=0}^{T-1}t^{{j\over T}}{\bf C}[t,t^{-1}]\otimes A_{j}.
\end{eqnarray}
We define a symmetric bilinear form on $\tilde{A}[\sigma]$ as follows:
\begin{eqnarray}
\langle t^{m+{j\over T}}\otimes a, t^{n+{k\over T}}\otimes b\rangle
=B(a,b)\delta_{m+n+1+{j+k\over T},0}
\end{eqnarray}
for $a\in A_{j}, b\in A_{k}, m,n\in \Bbb{Z}$. Let
$C(\tilde{A}[\sigma])$ be the Clifford algebra of $\tilde{A}[\sigma]$
with respect to this defined symmetric bilinear form.

Let $M$ be any restricted $C(\tilde{A}[\sigma])$-module, i.e., for any
$u\in M$, $(t^{n+{j\over T}}\otimes A_{j})u=0$ for sufficiently large
$n$. Then $\{
u(z)|u\in A\}$ is a local subspace of $F(M,T)$. Let $V$ be the vertex
superalgebra generated by $\{u(z)|u\in A\}$. Then $M$ is a faithful
$\sigma$-twisted $V$-module. By repeating our technique used in the first
subsection, we see that $V$ is a module for $C(\tilde{A})$.
Since $V$ contains $I(z)$, $V$ is a lowest weight $C(\tilde{A})$-module
so that $V=\Bbb{F}^{n}$. Thus $M$ is a $\sigma$-twisted
${\bf F}^{n}$-module.

Next, we follow [DM1] or [FFR] to find the irreducible
$C(\tilde{A}[\sigma])$-modules.
For $u\in A_{j}, n\in \Bbb{Z}$, we define
\begin{eqnarray}
\deg (t^{n+{j\over T}}\otimes u)=-n-{1\over 2}-{j\over T}.
\end{eqnarray}
Then $\tilde{A}[\sigma]$ becomes a ${1\over T}\Bbb{Z}$-graded (resp.
${1\over 2T}\Bbb{Z}$-graded) space if $T$ is even (resp. odd). Let
$\tilde{A}[\sigma]_{\pm}$ be the sum of homogeneous subspaces of positive or
negative degrees and let
$\tilde{A}[\sigma]_{0}$ be the degree-zero homogeneous subspace.
Then $\tilde{A}[\sigma]_{0}=0$ if $T$ is odd,
$\tilde{A}[\sigma]_{0}=t^{-{1\over 2}}\otimes A_{{T\over 2}}$ if $T$
is even. Then
$\tilde{A}[\sigma]_{\pm}$ are isotropic subspaces.
Let $\Lambda
(\tilde{A}[\sigma]_{+})$ be the exterior algebra of $\tilde{A}[\sigma]_{+}$.
If $T$ is odd, there is no nonzero element of degree-zero in the twisted
Clifford algebra.
Then it follows from the Stone-Von Neumann theorem that $\Lambda
(\tilde{A}[\sigma]_{+})$ is the unique irreducible
$C(\tilde{A}[\sigma])$-module and that if $M$ is a
$C(\tilde{A}[\sigma])$-module such that for any $u\in M$, there is a
positive integer $N$ such that $\tilde{A}[\sigma]_{-}^{N}u=0$, then
$M$ is a direct sum of some copies of $\Lambda
(\tilde{A}[\sigma]_{+})$.

If $T=2k$ is even, the  situation is a little bit complicated.
If $\dim A_{k}=2m$ is even, we can decompose
$A_{k}=(A_{k})_{+}\oplus (A_{k})_{-}$ where $(A_{k})_{\pm}$ are
isotropic subspaces of dimension $m$. Then
it follows from the Stone-Von Neumann theorem that $\Lambda
(\tilde{A}[\sigma]_{+}\oplus
(A_{k})_{+})$ is the unique irreducible
$C(\tilde{A}[\sigma])$-module up to equivalence, so that
$\Lambda (\tilde{A}[\sigma]_{+}\oplus
(A_{k})_{+})$ is the unique irreducible $\sigma$-twisted $\Bbb{
F}^{n}$-module up to equivalence.

If $\dim A_{k}=2m+1$ ($m\in \Bbb{Z}_{+}$) is odd, let $e\in A_{k}$
such that $\langle
e,e\rangle =2$. Then we have:
 \begin{eqnarray}
A_{k}={\bf C}e\oplus (A_{k})_{+}\oplus (A_{k})_{-}
\end{eqnarray}
where $(A_{k})_{\pm}$ are isotropic subspaces of dimension $m$ such
that $\langle e,(A_{k})_{\pm}\rangle=0$. Let
\begin{eqnarray}
\Lambda
(\tilde{A}[\sigma]_{+}\oplus (A_{k})_{+})=\Lambda^{even}
(\tilde{A}[\sigma]_{+}\oplus (A_{k})_{+})\oplus \Lambda^{odd}
(\tilde{A}[\sigma]_{+}\oplus (A_{k})_{+}).
\end{eqnarray}
Let $e$ act as $(\pm 1, \mp 1)$. Then we obtain
two irreducible $C(\tilde{A}[\sigma])$-modules, so that we get two
irreducible $\sigma$-twisted $\Bbb{F}^{n}$-modules. It is easy to prove that
these are the only irreducible modules, up to equivalence. Summarizing all the
arguments,
we have:

\begin{propo}
Any restricted $C(\tilde{A}[\sigma])$-module
$M$ is a weak $\sigma$-twisted $\Bbb{F}^{n}$-module. Furthermore, if
$T$ is odd,
$\Lambda (\tilde{A}[\sigma]_{+})$ is the unique irreducible
$\sigma$-twisted module for the vertex operator algebra $\Bbb{F}^{n}$;
if $T$ is even and $\dim A_{{T\over 2}}$ is even,
$\Lambda (\tilde{A}[\sigma]_{+}\oplus
(A_{k})_{+})$ is the unique irreducible $\sigma$-twisted $\Bbb{
F}^{n}$-module up to equivalence; if $T$ is even and $\dim A_{{T\over
2}}$ is odd, $\Bbb{F}^{n}$ has two irreducible $\sigma$-twisted
modules.
\end{propo}

\begin{propo}
Let $M$ be any $\sigma$-twisted weak $\Bbb{
F}^{n}$-module such that for any $u\in M$, there is a positive integer
$N$ such that $\tilde{A}[\sigma]_{-}^{N}u=0$. Then $M$ is a direct sum
of irreducible $\sigma$-twisted $\Bbb{F}^{n}$-modules $\Lambda
(\tilde{A}[\sigma]_{+})$.
\end{propo}

{\bf Case 2.} Let ${\bf g}$ be a finite-dimensional simple Lie algebra
with a fixed
Cartan subalgebra ${\bf h}$. Let $\Delta$ be the set of all roots, let
$\Pi $ be a set of simple roots, and let $\theta$ be the
highest root. Let $\langle\cdot,\cdot\rangle$ be the
normalized Killing form such that $\langle\theta,\theta\rangle=2$.
Then the Killing form is $2\Omega \langle\cdot,\cdot\rangle$, where
$\Omega$ is the dual Coxeter number of ${\bf g}$.
Let $\sigma$ be an automorphism of ${\bf g}$ such that $\sigma
^{T}=Id$. Then $\sigma$ preserves the invariant bilinear form.
Denote by $M_{({\bf g},\sigma)}(\ell,\lambda)$ the Verma module
for $\tilde{{\bf g}}[\sigma]$ with lowest weight $\lambda$ of level $\ell$.
Denote the irreducible quotient module by $L_{({\bf
g},\sigma)}(\ell,\lambda)$.

\begin{corol}
 Let $\sigma$ be any automorphism of order
$T$ of $M_{{\bf g}}(\ell,\Bbb{C})$
or ${\bf g}$ and let $\ell$ be any complex number such that $\ell\ne
-\Omega$. Then
any restricted $\tilde{{\bf g}}[\sigma]$-module of level $\ell$ is a
weak $\sigma$-twisted $M_{{\bf g}}(\ell,\Bbb{C})$-module. In
particular, each irreducible $\tilde{{\bf g}}[\sigma]$-module $L_{({\bf
g},\mu)}(\ell,\lambda)$ is a $\sigma$-twisted $M_{{\bf
g}}(\ell,\Bbb{C})$-module.
\end{corol}

In the next section we will study the $\sigma$-twisted module theory
for the vertex operator algebra $L_{{\bf g}}(\ell,0)$.

\section{ The semisimple twisted module theory for $L_{{\bf g}}(\ell,0)$}
Let ${\bf g}$ be a finite-dimensional simple Lie
algebra with a fixed Cartan subalgebra ${\bf h}$ and let $\Delta$ be
the set of roots. Fix a set
$\Pi$ of simple roots and let $\theta$ be the highest root.
Let  $\langle \cdot,\cdot\rangle$ be the
normalized Killing form such that $\langle \theta,\theta\rangle=2$.

Let $\sigma$ be an automorphism of ${\bf g}$ such that $\sigma^{T}=1$.
Then we may consider $\sigma$ as an automorphism of the vertex operator
algebra $L_{{\bf g}}(\ell,0)$ for any $\ell\ne -\Omega$.
In this section, we first prove that there is a Dynkin diagram
automorphism $\mu$ of ${\bf g}$ such that the equivalence classes of
$\sigma$-twisted $L_{{\bf g}}(\ell,0)$-modules one-to-one correspond
to the equivalence classes of $\mu$-twisted $L_{{\bf
g}}(\ell,0)$-modules. Then we prove that if $\ell$ is a positive integer,
then any lower truncated ${1\over
T}\Bbb{Z}$-graded weak $\mu$-twisted $L_{{\bf g}}(\ell,0)$ module
 is completely reducible and that
 the set of equivalence classes of irreducible $\mu$-twisted $L_{{\bf
g}}(\ell,0)$-modules is just the set of equivalence classes of
standard $\tilde{{\bf g}}[\mu]$-modules of level $\ell$.

{\bf Throughout
this section, $\ell$ will be fixed a positive integer.}

Recall the following proposition from [DL], [Li] or [MP] about untwisted
representation theory.

\begin{propo}\label{Li1}
 Let $e$ be any root
vector of ${\bf g}$ with root $\alpha$. If $L_{{\bf g}}(\ell,\lambda)$
is an integrable $\tilde{{\bf g}}$-module, then
$e(z)^{t\ell +1}=0$ acting on $L_{{\bf g}}(\ell,\lambda)$ where
$t=1$ if $\alpha$ is a long root, $t=2$ if $\alpha$ is a short root and ${
\bf g}\ne {\bf g}_{2}$, $t=3$ if $\alpha$ is a short root and ${\bf g}={\bf
g}_{2}$.
\end{propo}

The following proposition was proved in [DL], [FZ] and [Li].

\begin{propo}\label{Li}
 Let $\displaystyle{M=\oplus_{n\in \Bbb{Z}}M(n)}$ be any lower truncated
$\Bbb{Z}$-graded weak module for $L_{{\bf g}}(\ell,0)$ as a vertex operator
algebra.
Then $M$ is a direct sum of standard
$\tilde{{\bf g}}$-modules of level $\ell$.
\end{propo}

Recall the following theorem of Kac from [K]:

\begin{propo}\label{Kac2}
 Let ${\bf g}$ be a finite-dimensional
simple Lie algebra, let ${\bf h}$ be a Cartan subalgebra and let
$\Pi=\{\alpha_{1},\cdots,\alpha_{n}\}$ be a set of simple roots. Let
$\sigma \in {\rm Aut}{\bf g}$ be such that $\sigma ^{T}=1$. Then
$\sigma$ is conjugate to an automorphism of ${\bf g}$ in the form
\begin{eqnarray}
\mu {\rm exp}\left({\rm ad}\left({2\pi \sqrt{-1}\over T}h\right)\right),\;\;
h\in {\bf h}_{\bar{0}},
\end{eqnarray}
where $\mu$ is a diagram automorphism preserving ${\bf h}$ and $\Pi
$, ${\bf h}_{\bar{0}}$ is the fixed-point set of $\mu$ in
${\bf h}$, and $\langle \alpha _{i},h\rangle\in \Bbb{Z}$ $ (i=1,\cdots,n)$.
\end{propo}

If $\sigma _{1}$ and $\sigma_{2}$ are conjugate automorphisms of a
vertex operator algebra $V$, then it is clear
that the equivalence classes of $\sigma_{1}$-twisted $V$-modules
one-to-one correspond to the equivalence classes of
$\sigma_{2}$-twisted $V$-modules (see also [DM1]). By Proposition
5.3, we only need to study the twisted module theory for the automorphisms
in form (5.1). Let $G({\bf g})$ be the subgroup of ${\rm
Aut}{\bf g}$ generated by
all Dynkin diagram automorphisms of ${\bf g}$. Then $G({\bf g})=\Bbb{
Z}_{k}$ where $k=1$ for $B_{n},C_{n},G_{n}$, $k=2$ for
$A_{n},D_{n}, E_{n}$ $(n\ne 4)$; $k=3$ for $D_{4}$ [K].

Let $V$ be a vertex (operator) superalgebra. A {\it derivation} ([B], [Lia])
of $V$ is an endomorphism $f$ of $V$ such that
\begin{eqnarray}
f(Y(u,z)v)=Y(f(u),z)v+Y(u,z)f(v)\;\;\;\mbox{for any }u,v\in V.
\end{eqnarray}
Let $f$ be a derivation of $V$. Then
\begin{eqnarray}
e^{z_{0}f}Y(u,z)v=Y(e^{z_{0}f}u,z)e^{z_{0}f}v\;\;\;\mbox{for any
}u,v\in V.
\end{eqnarray}
If a derivation $f$ of $V$ is locally finite on $V$, i.e., for any
$u\in M$, $\{f^{j}u|j\in \Bbb{Z}_{+}\}$ linearly spans a
finite-dimensional subspace of $M$, then $e^{f}$ is
an automorphism of $V$.

Let $a$ be an even element of $V$. Then $a_{0}$ (the zero mode of the
vertex operator $Y(a,z)$) is a derivation of $V$.
Let $V$ be a vertex operator superalgebra and let
$a$ be an even element of $V$ with weight $k\le 1$. Since ${\rm
wt}a_{0}=k-1\le 0$, by the
two grading-restriction assumptions on $V$, $a_{0}$ is  a locally finite
derivation on $V$. Then $e^{a_{0}}$ is an automorphism of $V$ except
that $e^{a_{0}}$ may change the Virasoro element $\omega$.
Furthermore, let $a$ be
a primary element. Then
\begin{eqnarray}
a_{0}\omega={\rm Res}_{z}Y(a,z)\omega={\rm Res}_{z}e^{zL(-1)}Y(\omega,-z)a
=(k-1)L(-1)a.
\end{eqnarray}
Therefore, if $a$ is an even weight-one  primary element of a vertex operator
superalgebra $V$, then  $e^{a_{0}}\omega=\omega$, so that $e^{a_{0}}$ is an
automorphism of $V$.
Furthermore, suppose that $a_{0}$ is semisimple on $V$. Then
$e^{a_{0}}$ is of finite-order if and only if there is a positive
integer $T$ such that $\displaystyle{Spec\; (a_{0})\subseteq \frac{2\pi
i}{T}\Bbb{Z}}$.

Let $V$ be a vertex operator superalgebra, let $\sigma$ be an automorphism
of order $S$ of $V$  and let $h\in V$ be an even element such that
\begin{eqnarray}
L(n)h=\delta_{n,0}h,\; \sigma
(h)=h,\;[h_{m},h_{n}]=0\;\;\mbox{for }m,n\in \Bbb{Z}_{+}.
\end{eqnarray}
Furthermore, we assume that $h(0)$ is semisimple on $V$ such that
$Spec\;{h(0)}\subseteq {1\over T}\Bbb{Z}$ for some positive
integer $T$. Then $\sigma_{h}=\exp (2\pi \sqrt{-1}h(0))$ is an
automorphism of $V$ such that $\sigma_{h}^{T}=1$. Since $\sigma(h)=h$, we get
$\sigma h(0)=h(0)\sigma$.
Let $(M,Y_{M})$ be a $\sigma$-twisted $V$-module. For any $a\in V$, we define
\begin{eqnarray}
\bar{Y}_{M}(a,z)=Y_{M}\left(z^{h(0)}\exp
\left(\sum_{n=1}^{\infty}{h(n)\over -n}(-z)^{-n}\right)a,z\right).
\end{eqnarray}

\begin{propo}
 $(M,\bar{Y}_{M}(\cdot,z))$ is a  weak $(\sigma
\sigma_{h})$-twisted $V$-module.
\end{propo}

{\bf Proof.} For any $a\in V$, it follows from the grading-restriction
assumptions on $V$ that $\Delta (z)a$ is a finite sum of
rational powers of $z$, so that $\bar{Y}_{M}(a,z)$ satisfies the
truncation condition. Next, we check the $L(-1)$-derivative property.
Set
\begin{eqnarray}\Delta (z)=z^{h(0)}\exp
\left(\sum_{n=1}^{\infty}{h(n)\over -n}(-z)^{-n}\right).
\end{eqnarray}
Notice that $\Delta (z)$ is invertible. By the definition of $\Delta
(z)$, we get
\begin{eqnarray}{d\over dz}\Delta
(z)=\left(-\sum_{n=0}^{\infty}h(n)(-z)^{-n-1}\right)\Delta (z).
\end{eqnarray}
Since $[L(-1), h(0)]=0$ and
\begin{eqnarray}
\lbrack L(-1), \sum_{n=1}^{\infty}{h(n)\over
-n}(-z)^{-n}\rbrack=\sum_{n=1}^{\infty}h(n-1)
(-z)^{-n}=\sum_{n=0}^{\infty}h(n)(-z)^{-n-1},
\end{eqnarray}
we get
\begin{eqnarray}
[L(-1), \Delta (z)]=\!\left(\sum_{n=1}^{\infty}h(n-1)
(-z)^{-n}\right)\!\Delta (z)=\!\left(\sum_{n=0}^{\infty}h(n)(-z)^{-n-1}\right)
\!\Delta (z).
\end{eqnarray}
Therefore
\begin{eqnarray}
{d\over dz}\bar{Y}_{M}(a,z)&=&{d\over dz}Y_{M}(\Delta (z)a,z)\nonumber\\
&=&Y_{M}\left({d\over dz}\Delta (z)a,z\right)+{\partial\over\partial
z_{0}}Y_{M}(\Delta (z)a,z_{0})|_{z_{0}=z}\nonumber\\
&=&Y_{M}\left({d\over dz}\Delta (z)a,z\right)+Y_{M}(L(-1)\Delta (z)a,z)
\nonumber\\
&=&Y_{M}\left(\left({d\over dz}\Delta (z)+[L(-1),\Delta (z)]\right)a,z\right)
+\bar{Y}_{M}(L(-1)a,z)\nonumber\\
&=&\bar{Y}_{M}(L(-1)a,z).
\end{eqnarray}
For the $(\sigma\sigma_{h})$-twisted Jacobi identity, let us assume
that $\sigma a=\varepsilon_{S}^{j}a$ and $h(0)a={k\over T}a$. Then
\begin{eqnarray}
&
&\!z_{0}^{-1}\delta\!\left(\!\frac{z_{1}-z_{2}}{z_{0}}\!\right)\!\bar{Y}_{M}(a,z_{1})
\bar{Y}_{M}(b,z_{2})-
\varepsilon_{a,b}z_{0}^{-1}\delta\!\left(\!\frac{z_{2}-z_{1}}{-z_{0}}
\!\right)\!\bar{Y}_{M}(b,z_{2})\bar{Y}_{M}(a,z_{1})\nonumber\\
&=&\!z_{0}^{-1}\delta\left(\frac{z_{1}-z_{2}}{z_{0}}\right)Y_{M}(\Delta
(z_{1})a,z_{1})Y_{M}(\Delta (z_{2})b,z_{2})\nonumber\\
& &\!-\varepsilon_{a,b}z_{0}^{-1}\delta\left(\frac{z_{2}-z_{1}}{-z_{0}}
\right)Y_{M}(\Delta (z_{2})b,z_{2})Y_{M}(\Delta
(z_{1})a,z_{1})\nonumber\\
&=&\!z_{2}^{-1}\delta\left(\frac{z_{1}-z_{0}}{z_{2}}\right)
\left(\frac{z_{1}-z_{0}}{z_{2}}\right)^{-{j\over S}}
Y_{M}(Y(\Delta (z_{1})a,z_{0})\Delta (z_{2})b,z_{2}).
\end{eqnarray}
Our desired Jacobi identity is:
\begin{eqnarray}
&
&\!z_{0}^{-1}\delta\!\left(\!\frac{z_{1}-z_{2}}{z_{0}}\!\right)\!\bar{Y}_{M}(a,z_{1})
\bar{Y}_{M}(b,z_{2})-
\varepsilon_{a,b}z_{0}^{-1}\delta\!\left(\!\frac{z_{2}-z_{1}}{-z_{0}}
\!\right)\!\bar{Y}_{M}(b,z_{2})\bar{Y}_{M}(a,z_{1})\nonumber\\
&=&\!z_{2}^{-1}\delta\left(\frac{z_{1}-z_{0}}{z_{2}}\right)
\left(\frac{z_{1}-z_{0}}{z_{2}}\right)^{-{j\over S}-{k\over T}}
\bar{Y}_{M}(Y(a,z_{0})b,z_{2})\nonumber\\
&=&\!z_{2}^{-1}\delta\left(\frac{z_{1}-z_{0}}{z_{2}}\right)
\left(\frac{z_{1}-z_{0}}{z_{2}}\right)^{-{j\over S}-{k\over T}}
Y_{M}(\Delta (z_{2})Y(a,z_{0})b,z_{2}).
\end{eqnarray}
Therefore, it is equivalent to prove:
\begin{eqnarray}
& &z_{2}^{-1}\delta\left(\frac{z_{1}-z_{0}}{z_{2}}\right)
\left(\frac{z_{1}-z_{0}}{z_{2}}\right)^{-{j\over S}-{k\over T}}
\Delta (z_{2})Y(a,z_{0})\nonumber\\
&=&z_{2}^{-1}\delta\left(\frac{z_{1}-z_{0}}{z_{2}}\right)
\left(\frac{z_{1}-z_{0}}{z_{2}}\right)^{-{j\over S}}
Y(\Delta (z_{1})a,z_{0})\Delta (z_{2}),
\end{eqnarray}
or
\begin{eqnarray}
& &z_{2}^{-1}\delta\left(\frac{z_{1}-z_{0}}{z_{2}}\right)
\left(\frac{z_{1}-z_{0}}{z_{2}}\right)^{-{j\over S}-{k\over T}}
\Delta (z_{2})Y(a,z_{0})\Delta (z_{2})^{-1}\nonumber\\
&=&z_{2}^{-1}\delta\left(\frac{z_{1}-z_{0}}{z_{2}}\right)
\left(\frac{z_{1}-z_{0}}{z_{2}}\right)^{-{j\over S}}
Y(\Delta (z_{1})a,z_{0}).
\end{eqnarray}
{}From the commutator formula, we have
\begin{eqnarray}
[h(n), Y(a,z_{0})]=\sum_{i=0}^{\infty}\left(\begin{array}{c}n\\i\end{array}
\right)z_{0}^{n-i}Y(h(i)a,z_{0}).
\end{eqnarray}
Then
\begin{eqnarray}
& &\lbrack \sum_{n=1}^{\infty}\frac{h(n)}{-n}(-z_{2})^{-n}, Y(a,z_{0})\rbrack
\nonumber\\
&=&-\sum_{n=1}^{\infty}\sum_{i=0}^{\infty}{1\over n}\left(\begin{array}
{c}n\\i\end{array}\right)(-z_{2})^{-n}z_{0}^{n-i}Y(h(i)a,z_{0})\nonumber\\
&=&-\sum_{n=1}^{\infty}{1\over n}(-z_{2})^{-n}z_{0}^{n}Y(h(0)a,z_{0})
\nonumber\\
& &-\sum_{i=1}^{\infty}\sum_{n=1}^{\infty}{1\over i}\left(\begin{array}
{c}n-1\\i-1\end{array}\right)(-z_{2})^{-n}z_{0}^{n-i}Y(h(i)a,z_{0})\nonumber\\
&=&\log \left(1+{z_{0}\over z_{2}}\right)^{-{k\over T}}Y(a,z_{0})
+\sum_{i=1}^{\infty}\frac{(-z_{2}-z_{0})^{-i}}{-i}Y(h(i)a,z_{0}).
\end{eqnarray}
Thus
\begin{eqnarray}
& &\exp \left(\sum_{n=1}^{\infty}\frac{h(n)}{-n}(-z_{2})^{-n}\right)
Y(a,z_{0})\exp \left(-\sum_{n=1}^{\infty}\frac{h(n)}{-n}(-z_{2})^{-n}\right)
\nonumber\\
&=&Y\left(\left(\exp \log \left(1+{z_{0}\over z_{2}}\right)^{{k\over T}}\exp
\sum_{i=1}^{\infty}\frac{(-z_{2}-z_{0})^{-i}}{-i}h(i)\right)a,z_{0}\right)
\nonumber\\
&=&\left(1+{z_{0}\over z_{2}}\right)^{{k\over T}}Y\left((z_{2}+z_{0})^{-h(0)}
\Delta (z_{2}+z_{0})a,z_{0}\right)\nonumber\\
&=&z_{2}^{-{k\over T}}Y\left(\Delta (z_{2}+z_{0})a,z_{0}\right).
\end{eqnarray}
Since $[h(0),Y(a,z_{0})]=Y(h(0)a,z_{0})$, we get:
\begin{eqnarray}
z^{h(0)}Y(a,z_{0})z^{-h(0)}=Y(z^{h(0)}a,z_{0}).
\end{eqnarray}
Therefore
\begin{eqnarray}
\Delta(z_{2})Y(a,z_{0})\Delta (z_{2})^{-1}
&=&z_{2}^{-{k\over T}}Y\left(z_{2}^{h(0)}
\Delta (z_{2}+z_{0})a,z_{0}\right)\nonumber\\
&=&Y\left(\Delta (z_{2}+z_{0})a,z_{0}\right).
\end{eqnarray}
Then
\begin{eqnarray}
& &z_{2}^{-1}\delta\left(\frac{z_{1}-z_{0}}{z_{2}}\right)
\left(\frac{z_{1}-z_{0}}{z_{2}}\right)^{-{j\over S}-{k\over T}}
\Delta (z_{2})Y(a,z_{0})\Delta (z_{2})^{-1}\nonumber\\
&=&z_{1}^{-1}\delta\left(\frac{z_{2}+z_{0}}{z_{1}}\right)
\left(\frac{z_{2}+z_{0}}{z_{1}}\right)^{{j\over S}+{k\over T}}
Y\left(\Delta (z_{2}+z_{0})a,z_{0}\right)\nonumber\\
&=&z_{1}^{-1}\delta\left(\frac{z_{2}+z_{0}}{z_{1}}\right)
\left(\frac{z_{2}+z_{0}}{z_{1}}\right)^{{j\over S}}
Y\left(\Delta (z_{1})a,z_{0}\right).
\end{eqnarray}
Therefore, the $(\sigma\sigma_{h})$-twisted Jacobi identity holds. Thus
$(M,\bar{Y})$ is a weak $(\sigma\sigma_{h})$-twisted $V$-module.$\;\;\;\;\Box$

\begin{rema}
 Since $\Delta (\!z\!)$ is invertible, $(M,\bar{Y})$ is irreducible
if $(M,Y)$ is irreducible. In particular, if $\sigma
\!=\!Id_{V}$, any $\sigma_{h}$-twisted $V$-module can be constructed
{}from a $V$-module.
Proposition 5.4 also gives an isomorphism between the
twisted affine Lie algebras $\tilde{{\bf g}}[\mu]$ and $\tilde{{\bf
g}}[\mu\sigma_{h}]$. Although Proposition 5.4 is only about an inner
automorphism, we believe that it has more interesting applications.
We will study this subject in a coming paper.
\end{rema}

{}From Propositions 5.3 and 5.4, it is enough for us  to study the
$\mu$-twisted module
theory for $L_{{\bf g}}(\ell,0)$. Notice that ${\bf g}_{0}$ is a
simple subalgebra of ${\bf g}$ [K]. For any linear
functional $\lambda$ on ${\bf h}_{0}$, denote by $L_{({\bf
g},\mu)}(\ell,\lambda)$ the irreducible lowest weight $\tilde{{\bf
g}}[\mu]$-module of level $\ell$ with lowest weight $\lambda$. Recall
{}from [K] that $\theta_{0}$ is the highest weight of ${\bf g}_{1}$ with
respect to the Cartan subalgebra of ${\bf g}_{0}$ and that
$e_{\theta_{0}}$ is a root vector.

\begin{propo}
 An irreducible $\tilde{{\bf
g}}[\mu]$-module $L_{({\bf g},\mu)}(\ell,\lambda)$ is a
$\mu$-twisted module for the vertex operator algebra $L_{{\bf
g}}(\ell,0)$ if and only if it is a standard $\tilde{{\bf
g}}[\mu]$-module.
 \end{propo}

{\bf Proof}. For any $a\in {\bf g}_{j}$, by considering
$a(z)=\sum_{n\in \Bbb{Z}}(t^{n+{j\over T}}\otimes a)z^{-n-1-{j\over
T}}$ as an element of $F(L_{{\bf g},\mu}(\ell,\lambda),T)$, we obtain
a local subspace $\{a(z)|a\in {\bf g}\}$ of $F(L_{({\bf
g},\mu)}(\ell,\lambda),T)$.
Let $V$ be the vertex  algebra generated by
$\{a(z)|a\in {\bf g}\}$. Being a quotient vertex algebra of the vertex
operator algebra $M_{{\bf g}}(\ell,\Bbb{C})$, $V$ is a vertex operator
algebra. Then $L_{({\bf g},\mu)}(\ell,\lambda)$ is a
$\mu$-twisted $L_{{\bf g}}(\ell,0)$-module if and only if
 $V=L_{{\bf g}}(\ell,0)$. It follows from [Li] that $V=L_{{\bf
g}}(\ell,0)$ if and only if $Y_{V}(e_{\theta},z)^{\ell+1}=0$ on $V$. Since
$L_{({\bf g},\mu)}(\ell,\lambda)$ is a faithful $\mu$-twisted $V$-module,
it follows from Proposition 2.12 that
$Y_{V}(e_{\theta},z)^{\ell+1}=0$ on $V$ if and only if
$Y_{M}(e_{\theta},z)^{\ell+1}=0$ for $M=L_{({\bf g},\mu)}(\ell,\lambda)$.
Then we only need to prove that
$Y_{M}(e_{\theta},z)^{\ell+1}=0$ on $L_{({\bf
g},\mu)}(\ell,\lambda)$ if and only if $L_{({\bf
g},\mu)}(\ell,\lambda)$ is an integrable $\tilde{{\bf g}}[\mu]$-module.

If $L_{({\bf g},\mu)}(\ell,\lambda)$ is an integrable $\tilde{{\bf
g}}[\mu]$-module, then it is an integrable $\tilde{{\bf
g}}^{\theta}$-module. By Proposition 5.1, we have:
$Y_{M}(e_{\theta},z)^{\ell+1}=0$ for $M=L_{({\bf
g},\mu)}(\ell,\lambda)$. Then $L_{({\bf g},\mu)}(\ell,\lambda)$ is a
$\mu$-twisted module for the vertex operator algebra $L_{{\bf g}}(\ell,0)$.
Conversely, suppose that $L_{({\bf g},\mu)}(\ell,\lambda)$ is a
$\mu$-twisted $L_{{\bf g}}(\ell,0)$-module. From the untwisted theory [Li]
$L_{({\bf g},\mu)}(\ell,\lambda)$ is an integrable $\tilde{{\bf
g}_{0}}$-module.
By Proposition 5.1, $Y(e_{\theta_{0}},z)$ is nilpotent on $L_{{\bf
g}}(\ell,0)$. By Proposition 2.12, $Y(e_{\theta_{0}},z)$ is nilpotent
on $L_{({\bf g},\mu)}(\ell,\lambda)$. Therefore,
$L_{({\bf g},\mu)}(\ell,\lambda)$ is
an integrable $\tilde{{\bf g}}[\mu]$-module. $\;\;\;\;\Box$

Let us recall a complete reducibility theorem of Kac's (Theorem 10.7 [K]).

\begin{propo}\label{Kac1}
Let $g(A)$ be a Kac-Moody algebra
associated to a symmetrilizable generalized Cartan matrix $A$ of rank
$n$, let $e_{i}, f_{i}, h_{i}$ $(i=1,\cdots,n)$ be the Chevalley
generators for the derived subalgebra $g'(A)$ and let
$g(A)=g(A)_{+}\oplus H\oplus g(A)_{-}$ be the
triangular decomposition.
Let $M$ be a $g'(A)$-module such that for any $u\in M$
there is a positive integer $k$ such that
\begin{eqnarray}
g(A)_{+}^{k}u=0,\;f_{i}^{k}u=0\;\;\;\mbox{{\it for
 }}i=1,\cdots,n.
\end{eqnarray}
Then $M$ is a direct sum of irreducible highest weight integrable
$g'(A)$-modules.
\end{propo}

For our purpose, we are only concerned about (both
twisted and untwisted) affine Lie algebras. Let ${\bf g}$ be a
finite-dimensional simple
Lie algebra, let $\mu$ be a Dynkin diagram automorphism
of order $T$ of ${\bf g}$ and let $\tilde{{\bf g}}[\mu]$ be the
corresponding twisted affine Lie algebra as before. Then we have the
following modified complete reducibility theorem.

\begin{propo}
Let $\tilde{{\bf g}}[\mu]$ be the twisted affine
Lie algebra as above and let $M$ be a $\tilde{{\bf g}}[\mu]$-module
such that for any $u\in M$, there are positive integers $k$ and $r$ such that
\begin{eqnarray}
\tilde{{\bf g}}[\mu]_{+}^{k}u=0,\; e_{\alpha}(z)^{r}M=0
\end{eqnarray}
where $\alpha$ is a root of ${\bf g}_{0}$, $e_{\alpha}$ is a root
vector of root $\alpha$. Then $M$
is a direct sum of irreducible highest weight integrable $\tilde{{\bf
g}}[\mu]$-modules (of level less than $r$).
\end{propo}

{\bf Proof.} For such a $\tilde{{\bf g}}[\mu]$-module $M$, we set
\begin{eqnarray}
\Omega (M)=\{u\in M|\tilde{{\bf g}}[\mu]_{+}u=0\}.
\end{eqnarray}
Then $\Omega (M)\ne 0$ if $M\ne 0$. It is also clear that $\Omega (M)$
is a ${\bf g}_{0}$-module. Since $e_{\alpha}(z)^{r}\Omega (M)=0$, by
considering the coefficient of $z^{-r}$ we obtain
$e_{\alpha}(0)^{r}\Omega (M)=0$. It follows from Proposition 5.1.5 [Li]
that $\Omega (M)$ is a direct sum of finite-dimensional irreducible
${\bf g}_{0}$-modules. Let $u$ be a highest weight vector of $\Omega
(M)$ as a ${\bf g}_{0}$-module. Then by considering the constant term
of $e_{\alpha}(z)^{r}u=0$ for a positive root $\alpha$ we obtain
$e_{\alpha}(-1)^{r}u=0$. Thus the submodule $M^{o}$ of $M$ generated by
$\tilde{{\bf g}}[\mu]$ from $\Omega (M)$ satisfies the conditions
assumed in Proposition 5.7, so that $M^{o}$ is a direct sum of
irreducible highest weight integrable $\tilde{{\bf g}}[\mu]$-modules.

Next we shall prove that $M=M^{o}$. Suppose $M\ne M^{o}$. Then $\Omega
(M/M^{o})\ne 0$, so that $(M/M^{o})$ has a (nonzero) irreducible
highest weight integrable submodule $M^{1}/M^{o}$. It is easy to see
that $M^{1}$ satisfies the conditions assumed in Proposition 5.7. Thus
$M^{1}$ is completely reducible, so that $M^{1}=M^{o}\oplus W$ where
$W$ is a standard module. But by definition the highest weight vectors of
$W$ are contained in $M^{o}$. This is a contradiction.$\;\;\;\;\Box$

\begin{theo}
 Let $M$ be a   $\mu$-twisted weak
$L_{{\bf g}}(\ell,0)$-module such that for any $u\in M$ there is a positive
integer $k$ such that $\tilde{{\bf g}}[\mu]_{+}^{k}u=0$. Then $M$ is a
direct sum of irreducible $\mu$-twisted $L_{{\bf g}}(\ell,0)$-modules,
which are standard $\tilde{{\bf g}}[\mu]$-modules of level $\ell$.
\end{theo}

{\bf Proof}. It follows directly from Propositions 5.1, 5.6 and 5.8.
$\;\;\;\;\Box$

Since any lower truncated ${1\over T}\Bbb{Z}$-graded
weak  $\sigma$-twisted $L_{{\bf g}}(\ell,0)$-module $M$
satisfies the conditions of Theorem 5.9, we have:

\begin{corol}
Any lower truncated ${1\over T}\Bbb{Z}$-graded
weak  $\sigma$-twisted $L_{{\bf g}}(\ell,0)$-module $M$
is a direct sum of standard $\tilde{{\bf
g}}[\sigma]$-modules of level $\ell$.
\end{corol}

By Proposition 5.4 and Remark 5.5 we have:
\begin{theo}
Let $\sigma$ be any automorphism of finite order
of ${\bf g}$ and let $M$ be any weak $\sigma$-twisted  $L_{{\bf
g}}(\ell,0)$-module such that for any $u\in M$ there is a positive
integer $k$ such that $\tilde{{\bf g}}[\sigma]_{+}^{k}u=0$. Then $M$
is completely reducible and there are only finite many irreducible
$\sigma$-twisted $L_{{\bf g}}(\ell,0)$-modules up to equivalence.
\end{theo}

For the rest of this section, we apply Proposition 5.4 to the study of
certain twisted $\Bbb{F}^{2n}$-modules. Following [DM1], let
$A=A_{+}\oplus A_{-}$ and define
$G=GL(A^{+})$. For any $\sigma\in G$, we define an action of $\sigma$
on $A^{-}$ by:
\begin{eqnarray}
\langle u,\sigma v\rangle=\langle \sigma^{-1}u,v\rangle
\;\;\;\mbox{ for any }u\in A^{+},v\in A^{-}.
\end{eqnarray}
Then $G$ acts on $A$ such that $G$ preserves the bilinear form on $A$.
For $\sigma\in G$ of order $T$, let
$\{a_{1,\sigma},\cdots,a_{n,\sigma}\}$ be a basis of $A^{+}$ such that
\begin{eqnarray}
\sigma a_{i,\sigma}=\varepsilon^{n_{i,\sigma}}a_{i,\sigma}
\end{eqnarray}
where $0\le n_{i,\sigma}<T$ for each $1\le i\le n$. Let $\{
a_{1,\sigma}^{*},\cdots,a^{*}_{n,\sigma}\}$ be a dual basis of $A^{-}$
so that $\langle a_{i,\sigma},a_{j,\sigma}^{*}\rangle=\delta_{i,j}$. Then
\begin{eqnarray}
\sigma a^{*}_{i,\sigma}=\varepsilon^{-n_{i,\sigma}}a^{*}_{i,\sigma}.
\end{eqnarray}

Set
\begin{eqnarray}
h_{\sigma}={1\over
T}\sum_{j=1}^{n}n_{j,\sigma}(a_{j,\sigma})_{-1}a^{*}_{j,\sigma}.
\end{eqnarray}
Then $h_{\sigma}$ is an even weight-one element of $\Bbb{F}^{2n}$. For
any $a,b,c\in A$, we have:
\begin{eqnarray}
& &\!\!(a_{-1}b)_{0}c\nonumber\\
&=&\!\!{\rm Res}_{z_{1}}{\rm Res}_{z_{2}}\left((z_{1}-z_{2})^{-1}Y(a,z_{1})
Y(b,z_{2})+(-z_{2}+z_{1})^{-1}Y(b,z_{2})Y(a,z_{1})\right)c\nonumber\\
&=&\!\!\sum_{r=0}^{\infty}\left(a_{-r-1}b_{r}c-b_{-r-1}a_{r}c\right)\nonumber\\
&=&\!\!a_{-1}b_{0}c-b_{-1}a_{0}c\nonumber\\
&=&\!\!\langle b,c\rangle a-\langle a,c\rangle b.
\end{eqnarray}
Then
\begin{eqnarray}
h_{\sigma}(0)a_{i,\sigma}={1\over T}n_{i,\sigma}a_{i,\sigma},\;
h_{\sigma}(0)a^{*}_{i,\sigma}=-{1\over T}n_{i,\sigma}a^{*}_{i,\sigma}.
\end{eqnarray}
Since $A$ generates the vertex operator superalgebra $\Bbb{F}^{2n}$,
we have:
\begin{eqnarray}
\sigma =\exp\left(2\pi \sqrt{-1}h_{\sigma}(0)\right).
\end{eqnarray}
(See [DM1].) Then as a corollary of Proposition 5.4, we have the following
 proposition of Dong and Mason [DM1]:

\begin{propo}
For any $\sigma\in G$,  up to equivalence
there is only one irreducible $\sigma$-twisted $\Bbb{F}^{2n}$-module.
Furthermore,
it can be constructed by using Proposition 5.4 from the adjoint $\Bbb{
F}^{2n}$-module.
\end{propo}

\end{document}